\documentclass{amsart}
\usepackage{amsmath, amssymb, amsfonts, amsthm, mathrsfs, hyperref, enumitem, graphicx, xcolor}
\usepackage[square,numbers]{natbib}
\newtheorem{theorem}{Theorem}
\usepackage{epstopdf}
\newtheorem{lemma}{Lemma}
\newtheorem{definition}{Definition}
\newcommand{\bZ}{\mathbb{Z}}
\newcommand{\bR}{\mathbb{R}}
\newcommand{\bN}{\mathbb{N}}

\begin{document}
	
	\title[Short Simple Geodesic Loops]{Short Simple Geodesic Loops on a 2-Sphere}
\author[Beach]{Isabel Beach}
\address[Beach]{Department of Mathematics,  University of Toronto, Toronto, Canada}
\email{isabel.beach@mail.utoronto.ca}
	
	\begin{abstract}
		The classic Lusternik--Schnirelmann theorem states that there are three distinct simple periodic geodesics on any Riemannian 2-sphere $M$. It has been proven by Y. Liokumovich, A. Nabutovsky and R. Rotman that the shortest three such curves have lengths bounded in terms of the diameter $d$ of $M$. We show that at any point $p$ on $M$ there exist at least two distinct simple geodesic loops (geodesic segments that start and end at $p$) whose lengths are respectively bounded by $8d$ and $14d$.
	\end{abstract}
	
	\maketitle
	
\section{Introduction}
\sloppy
A geodesic loop based at a point $p$ in a Riemannian manifold $M$ is a geodesic $\gamma:[0,1]\to M$ with $\gamma(0)=\gamma(1)=p$. At each point $p$ on a Riemannian 2-sphere $M$, there are infinitely many distinct geodesic loops based at $p$, as shown by J. P. Serre in \cite{serre1951}. Length bounds for geodesic loops have been studied by many authors, including S. Sabourau in \cite{sabourau2004} and F. Balacheff, H. Parlier and Sabourau in \cite{balacheff_2010_loops}. In \cite{rotman_2009, rotman_2011_loops, rotman_2008_loops}, A. Nabutovsky and R. Rotman produced length bounds for these geodesic loops in terms of the geometric properties of $M$, culminating in the result that on a Riemannian manifold $M$ with diameter $d$ diffeomorphic to $S^2$, there are $k$ distinct geodesic loops based at any $p\in M$ of length at most $20kd$ for every positive integer $k$ \cite{rotman_2011_loops}. Building on this work, H. Y. Cheng improved the bound to $6kd$ in general and to $5kd$ on a generic set of metrics on $S^2$ when $k$ is odd \cite{cheng2022}. Rotman also proved in \cite{rotman_2008_loops} that the shortest geodesic loop at any point on a closed Riemannian manifold of dimension $n$ has length at most $2nd$. Thus the shortest two non-trivial geodesic loops at any point on a Riemannian 2-sphere are known to have lengths at most $4d$ and $12d$. 
\par 
Similarly, a closed (or periodic) geodesic is a geodesic curve $\gamma:\bR\to M$ that is also periodic, and hence has an image that can be parameterized by a geodesic loop. In this context, Rotman and Nabutovsky proved a quantitative version of the Lusternik--Schnirelmann theorem, showing that on a Riemannian 2-sphere of diameter $d$ there are always at least two simple closed geodesics of respective lengths at most $5d$ and $10d$ \cite{rotman_ls_2011}. With Y. Liokumovich, they upgraded this result to prove the existence of a third simple closed geodesic of length at most $20d$ \cite{rotman_ls_2017}. Note that because these results deal with closed geodesics, they do not produce curves with a specified base point.
\par
In this paper, we are interested in the existence of simple geodesic loops of bounded length. We will prove the following theorem.
\begin{theorem}
	Let $p$ be a point in an analytic Riemannian 2-sphere $M$ of diameter $d$. Then there are at least two simple geodesic loops based at $p$ of respective lengths at most $8d$ and $14d$.
	\label{theorem:main}
\end{theorem}
\noindent 
This differs from previous results regarding geodesic loops because it concerns simple curves. Consequently, the two loops obtained in the proof of our theorem will have distinct images. This is in contrast to the work of Cheng, Nabutovsky and Rotman, which does not preclude the possibility that the $k$ short geodesic loops are given by iterating the same closed geodesic $k$ times. This result also differs from the quantitative Lusternik--Schnirelmann theorem in that it produces simple loops based at any point in $M$.
\par
Our proof is built on the Lusternik--Schnirelmann proof of the existence of three simple geodesics on the 2-sphere and the quantitative techniques used by Nabutovsky and Rotman. We follow the versions of the Lusternik--Schnirelmann proof given by J. Hass and P. Scott in \cite{hass1994} and W. Klingenberg in \cite{klingenberg_1978}. In order to produce three simple closed geodesics, the standard technique is to use representatives of three non-trivial $\bZ_2$ homology classes in the space of unparameterized, unoriented simple closed curves in $M$. In the cited proofs, one then iteratively applies a shortening process to the curves in the image of each representative. A subsequence of curves from each representative will then eventually converge to a non-trivial closed geodesic. Using the fact that these three homology classes are related by the cup product, it can be shown that either the three limiting geodesics obtained in this manner are distinct, or that there are infinitely many such geodesics.
\par 
Our paper is organized as follows. In Section \ref{sec:modified_disk_flow}, for any fixed $p\in M$ we develop a curve shortening procedure that depends on $p$ and is based on the disk flow of J. Hass and P. Scott \cite{hass1994}. We call our new process the modified disk flow at $p$. Instead of unparameterized simple closed curves, we will apply this process to what we call the $p$-admissible curves (see Definition \ref{def:admissible}). Roughly speaking, these are parameterized closed curves based at $p$ that can be approximated by simple curves. We will prove that the modified disk flow preserves the set of $p$-admissible curves (Lemma \ref{lemma:rays_non_transverse}) and does not increase the lengths of these curves (Lemma \ref{lemma:rays_are_shorter}). Importantly, we will show that the image of a $p$-admissible curve under this process converges (up to a subsequence) to a union of simple geodesic loops based at $p$ (Lemma \ref{lemma:modified_disk_flow}). 
\par
Next, in Section \ref{sec:short_meridional_slicings} we discuss how to construct a continuous map of $S^2$ into $M$ of non-zero degree that sends the meridian curves on $S^2$ to curves on $M$ that start at $p$ and end at some other point $q$. We call this map a short meridional slicing of $M$, and we will use these curves to form two homology cycles in the space of $p$-admissible curves. We first explain how to shorten a closed curve consisting of a pair of minimizing geodesics by applying the modified disk flow and exploiting the structure of the cut locus (Lemma \ref{lemma:digons_and_loops}). We then pick a point $q$ at locally maximal distance from $p$ and divide $M$ into disks bounded by minimizing geodesics connecting $p$ to $q$ using Berger's Lemma (see Lemma 4.1 in Chapter 13 Section 4 of \cite{docarmo}). Using our previous lemma, we contract the boundaries of these disks and turn the resulting homotopies into a short meridional slicing of $M$ (Lemma \ref{lemma:really_short_slicing}). Critically, the lengths of these curves will be bounded by the diameter of $M$. 
\par 
In Section \ref{sec:families}, we describe how to extend the modified disk flow to certain families of unparameterized, unoriented $p$-admissible curves. The modified disk flow on a single curve is defined with respect to a prescribed cover of $M$ by metric balls. A major obstacle to defining the modified disk flow on a family of curves is the existence of curves that are tangent to the boundary of a metric ball in this cover, because such curves can possibly give rise to discontinuities in the modified disk flow. First we show how to ensure that the set of parameters where such discontinuities occur is triangulable (Lemma \ref{lemma:analytic}). We then show how to remove these discontinuities (Lemma \ref{lemma:geo_polys_epsilon}). Lastly, we define the modified disk flow on a multi-parameter family and prove that the resulting family is in the same homotopy class and consists of curves of bounded length (Theorem \ref{theorem:flow}).
\par 
Finally, in Section  \ref{sec:counting_loops}, we prove the existence of a pair of short simple geodesic loops at $p$. We first show how to convert our short meridional slicing into a sweepout of $M$ by short $p$-admissible curves. After some minor modifications, we will then have two non-trivial cycles: a one-dimensional cycle formed by the slicing loops, and a two-dimensional cycle formed by pairs of slicing loops. We conclude by applying the modified disk flow to the images of the cycles (Lemma \ref{lemma:final_lemma}). Because these cycles are non-trivial, some subsequence of curves in each image converges to a union of short, non-trivial simple geodesic loops based at $p$. Following the Lusternik--Schnirelmann proof, we then show that there are at least two distinct such loops.

\section{Results}

\subsection{The Modified Disk Flow}
\label{sec:modified_disk_flow}
Our first step will be to create our curve shortening procedure. The original disk flow developed by Hass and Scott in \cite{hass1994} proceeds on a compact manifold $M$ as follows. Supposing we wish to shorten a closed curve $\gamma$, the disk flow produces a family of curves $\gamma_t$ for $t\geq 0$ with $\gamma_0=\gamma$. First, cover $M$ by a finite collection of totally normal metric balls $\{B_i\}_{i=1}^n$ (i.e., metric balls with the property that every point in their closure is connected by a unique minimizing geodesic that lies in the ball itself). In addition, the collection of metric balls with the same centres but half the radii should also cover $M$ and be in general position, so that no three balls intersect at a point. Note that in this paper we will consider all metric balls to be closed. We will consider the balls $B_i$ consecutively, starting at $B_1$. If $\gamma\cap B_1$ is empty, we do nothing. If instead $\gamma$ lies entirely in $B_1$, then it can be contracted to a point curve in $B_1$. Otherwise, consider each arc of $\gamma\cap B_1$ that connects two points on $\partial B_1$. Because $B_1$ is totally normal, there is a unique geodesic connecting the endpoints of each such arc. Moreover, this geodesic, excluding its endpoints, lies in $B_1$. We redefine $\gamma$ to be the curve obtained by replacing the arcs of $\gamma\cap B_1$ with these geodesics. We then repeat this process in each consecutive ball. The resulting curve is defined to be $\gamma_1$. We define $\gamma_t$ for $t\in(0,1)$ to be a suitable homotopy between $\gamma$ and $\gamma_1$ -- namely, a homotopy that does not increase the number of self-intersections. We then iterate this entire procedure to define $\gamma_t$ for $t\in[n,n+1]$, $n\in\mathbb{N}$, by repeating the above with $\gamma_n$ in place of $\gamma$.
\par 
When straightening in $B_i$, if the original arcs of $\gamma\cap B_i$ do not intersect then neither will the minimizing geodesic segments they are replaced with. This is because two minimizing geodesic segments cannot mutually intersect more than once in their interiors, and they cannot intersect non-transversely. Moreover, since our interpolating homotopy does not increase the number of self-intersections, the disk flow preserves simple curves. Thus it is possible to obtain a simple closed geodesic by starting with a simple curve.
\par 	
This flow can also be extended to multi-parameter families of curves. However, discontinuities in the disk flow can occur when one curve being shortened is tangent to the boundary of one of the balls in the cover and a nearby curve lies entirely within said ball. While this situation can be avoided for a single curve by slightly altering the cover, once we pass to a family of curves we cannot necessarily arrange the disks to avoid tangencies in every curve. Instead, each discontinuity is fixed by ``filling in the gap'' with an interpolating homotopy. We will discuss this process in detail in Section \ref{sec:families}.
\par 
Because we are interested in finding loops based at $p$, we want to modify the disk flow so that it fixes the base point of the curve we are shortening. Although it would be ideal to deal only with simple curves based at $p$, we will need to examine what we call $p$-admissible curves.
\begin{definition}[$p$-Admissible Curves]
	\label{def:admissible}
	Let $p\in M$. We call a curve $\gamma:[0,1]\to M$ $p$-admissible if the following conditions hold.
	\begin{enumerate}
		\item 
		$\gamma(0)=\gamma(1)=p$.
        \item $\gamma$ is piecewise geodesic and parameterized by arc length.
		\item 
		$\gamma$ has no self-intersections except perhaps at the point $p$.
		\item 
		There exists a homotopy of arbitrarily small width between $\gamma$ and a simple closed curve.
	\end{enumerate}
\end{definition}
\noindent In particular, such curves have no transverse self-intersections.
\par
\begin{figure}[ht]
	\centering
	\includegraphics[width=\textwidth]{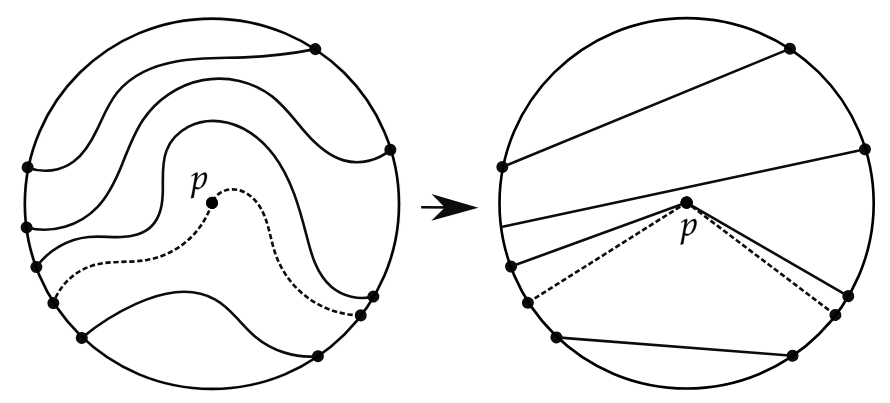}
	\caption{An example of the modified disk flow in $B_1$. The dashed arc contains the point corresponding to $t=0$.}
	\label{fig:rays}
\end{figure}

We now describe how we will alter the disk flow to create what we call the modified disk flow at $p$. Note that we will only define this flow on $p$-admissible curves. We cover $M$ by balls $B_i$ as before. However, we additionally require that the first ball $B_1$ in our cover is centred at $p$, and that no other ball contains $p$. We will only make modifications to the flow in $B_1$. In order to fix $p=\gamma(0)$, we would like to replace the arc of $\gamma$ that contains $t=0$ with the two geodesics that connect $p$ to the arc's endpoints (we call such geodesics ``rays" because they emanate from the centre of the metric ball). However, unlike in the standard disk flow, this can create new transverse self-intersections since it is possible that another geodesic segment crosses both geodesic rays. Therefore we make the following compromise. First, the arc in $B_1$ containing $t=0$ is always replaced by a pair of geodesic rays. 
Secondly, if replacing an arc in $B_1$ by a minimizing geodesic would create transverse self-intersections with the pair of rays incident to $\gamma(0)$, we instead replace it with the pair of geodesic rays that connect its endpoints to $p$. All other arcs are replaced by genuine minimizing geodesics as in the original disk flow (see Figure \ref{fig:rays}).
\par
\begin{figure}
	\centering
	\includegraphics[width=\textwidth]{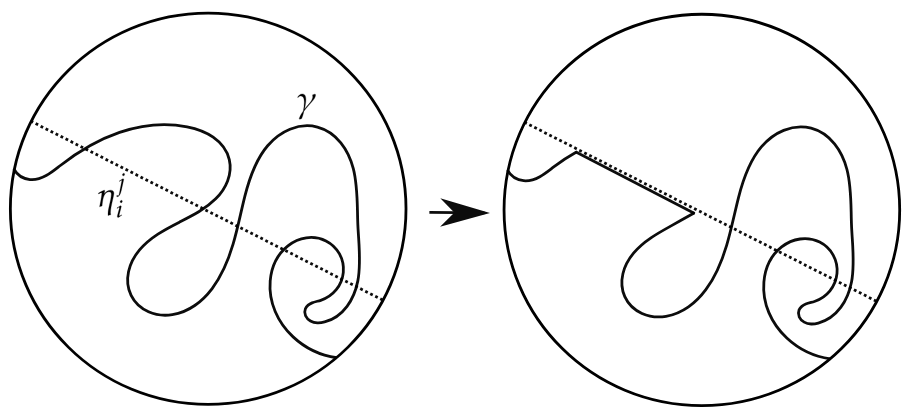}
	\caption{Eliminating a digon in order to homotope an arc to the corresponding geodesic segment. The dashed line is $\eta_i^j$ and the solid line is an arc of $\gamma$.}
	\label{fig:sliding_digons}
\end{figure}

 
The above determines $\gamma_1$, which we call the image of $\gamma$ under the modified disk flow. We define $\gamma_t$ for $t\in(0,1)$ to be a suitable homotopy between $\gamma=\gamma_0$ and $\gamma_1$. The homotopy $\gamma_t$ is in fact the same as in the original disk flow, despite the fact that we allow geodesic rays. The original description of the homotopy can be found in Lemma 1.6 of \cite{hass1994}. For completeness, we briefly reprove it here for $p$-admissible curves.
\begin{lemma}
	\label{lemma:sliding_digons}
	A $p$-admissible curve can be homotoped to its image under the modified disk flow through $p$-admissible curves.
\end{lemma}
\begin{proof}
	For brevity, let the image of $\gamma$ under the modified disk flow at $p$ be denoted by $\eta$ (i.e, $\eta=\gamma_1$). It is enough to define a homotopy in $B_i$ between $\gamma\cap B_i$ and $\eta\cap B_i$ for each $i$. Enumerate the arcs of $\gamma\cap B_i$ and call them $\gamma_i^j$, and let the corresponding arcs of $\eta\cap B_i$ be called $\eta_i^j$. 
	To ensure our homotopy is through $p$-admissible curves, we want to homotope $\gamma_i^j$ to $\eta_i^j$ through arcs that do not intersect each other except perhaps non-transversely at $p$.
	\par
	If the region(s) bounded by $\gamma_i^j*-\eta_i^j$ contained no other arcs of $\gamma$, we could homotope $\gamma_i^j$ to $\eta_i^j$ without creating intersections. By $p$-admissibility, if any other $\gamma_i^k$ crosses $\eta_i^j$, $k\not=j$, necessarily they cross an even number of times. Thus the region bounded by $\gamma_i^j*-\eta_i^j$ is divided into regions bounded by an arc of $\eta_i^j$ and an arc of $\gamma$. Call these regions digons.  
	Of all the digons cobounded by arcs of $\eta_i^j$ (i.e., bounded by the union of an arc of $\eta_i^j$ and some other arc), consider one that does not contain any other digon cobounded by an arc of $\eta_i^j$. We can homotope the other bounding arc of this digon to an arc of $\eta_i^j$ without intersecting any other arc of $\gamma$ (see Figure \ref{fig:sliding_digons}), thus removing this digon. Then we slide the arc slightly so it lies outside the bounds of $\gamma_i^j*-\eta_i^j$ and does not intersect $\eta_i^j$. Repeat this process until all the digons cobounded by $\eta_i^j$ inside $\gamma_i^j*-\eta_i^j$ have been removed. We can then homotope $\gamma_i^j$ to $\eta_i^j$ without introducing intersection points. Repeat this process for every $j$.
\end{proof}    
\noindent
This defines $\gamma_t$ for $t\in[0,1]$. As before, we define $\gamma_t$ for $t\in[n,n+1]$, $n\in\mathbb{N}$ by iterating the modified disk flow procedure described above.
\par
\begin{figure}[ht]
	\centering
	\includegraphics[width=0.5\textwidth]{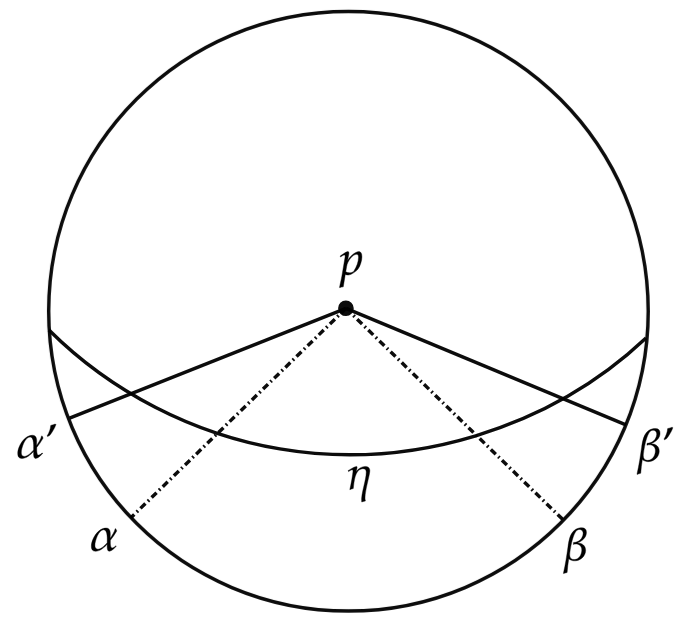}
	\caption{A minimizing geodesic $\eta$ that crosses the pair of rays $\alpha'$ and $\beta'$ must also cross the inner pair of rays $\alpha$ and $\beta$.}
	\label{fig:rays_non_transverse}
\end{figure}

 
Now that we have defined our procedure, we prove a lemma that clarifies some of its behaviour.
\begin{lemma}
	\label{lemma:rays_acute_angle}
	Let $\gamma$ be a $p$-admissible curve in $M$. Suppose the arc $\gamma_0$ of $\gamma\cap B_1$ containing $\gamma(0)$ intersects $\partial B_1$ in two points. Suppose there is another arc of $\gamma\cap B_1$, say $\gamma'$. Let $B'$ be the arc of $\partial B_1$ connecting the endpoints of $\gamma'$ that does not contain the endpoints of $\gamma_0$. The arc $B'$ subtends an angle $\theta$ in $\partial B_1$. Under the modified disk flow, $\gamma'$ will be replaced by a minimizing geodesic if and only if $\theta\leq\pi$.
\end{lemma}
\begin{proof}
	Note that such an arc $B'$ exists because $\gamma$ is $p$-admissible. 
	Denote the geodesic rays connecting $p$ to the endpoints of $\gamma_0$ by $\alpha$ and $\beta$. We need to determine when the minimizing geodesic $\eta$ connecting the endpoints of $\gamma'$ crosses $\alpha\cup\beta$. 
	Denote the geodesic rays connecting $p$ to the endpoints of $\gamma'$ by $\alpha'$ and $\beta'$. If $\theta\leq \pi$, then $\alpha'\cup\beta'\cup B'$ bounds a convex region, and hence the closure of this region contains $\eta$. On the other hand, this region does not contain $\alpha\cup\beta$, so $\eta$ will not intersect them (except perhaps non-transversely at $p$ if $\theta=\pi$). Thus $\gamma'$ will be replaced by $\eta$ under the modified disk flow. 
	Alternatively, if $\theta>\pi$, then the complement of the region bounded by $\alpha\cup\beta$ is convex, and hence the complement contains both $\eta$ and $\alpha\cup\beta$. Necessarily, $\eta$ must cross both $\alpha$ and $\beta$ transversely. Thus $\gamma'$ will be replaced by two rays through $p$ under the modified disk flow.
\end{proof}
The rest of this section is dedicated to proving certain properties of the modified disk flow. In order to successfully apply the Lusternik--Schnirelmann proof, we require the following.	
\begin{enumerate}
	\item 
	The modified disk flow must preserve $p$-admissibility.
	\item 
	Under repeated application of the modified disk flow, a curve must converge (up to a subsequence) to either the point curve $p$, to a union of at least two distinct simple geodesic loops at $p$, or to a single prime simple geodesic loop at $p$.
	\item 
	Under the modified disk flow, a curve's length is strictly decreased unless it is a union of geodesic loops at $p$.
	\item 
	The modified disk flow can be extended to one- and two-parameter families of $p$-admissible curves.
	\item 
	A family of $p$-admissible curves must be homotopic via $p$-admissible curves to its image under the modified disk flow.
\end{enumerate}
We now prove the first three properties. The last two properties will be proven in Section \ref{sec:families}.
We first prove that the modified disk flow preserves $p$-admissibility.
\begin{lemma}
	\label{lemma:rays_non_transverse}
	Let $\gamma$ be a $p$-admissible curve in $M$. Then the image of $\gamma$ under the modified disk flow is also $p$-admissible.
\end{lemma}
\begin{proof}
	Because we are replacing arcs with minimizing geodesics, the original disk flow does not create transverse intersections. This is shown in Theorem 1.8 of \cite{hass1994}. Therefore, because $p$-admissible curves have self-intersections only in $B_1$ and the only differences between the original and modified disk flows occur in $B_1$, we will only consider what happens in $B_1$. If $\gamma$ lies entirely in $\overline{B_1}$, then its image under the disk flow is the point $p$. Otherwise, the arc $\gamma_0$ of $\gamma\cap B_1$ containing $\gamma(0)$ intersects $\partial B_1$ in two distinct points. Denote the two geodesic rays that connect these two points to $p$ by $\alpha$ and $\beta$.
	\par 
	Apply the modified disk flow to $\gamma$ in $B_1$. Because $\gamma$ is $p$-admissible, if we obtain a transverse self-intersection after the modified disk flow it must be due to a minimizing geodesic crossing a pair of rays. By Lemma \ref{lemma:rays_acute_angle}, every ray pair obtained by applying the modified disk flow bounds a convex sector of $B_1$ that contains $\alpha$ and $\beta$. Therefore a minimizing geodesic that crosses any ray pair transversely must also cross $\alpha$ and $\beta$ transversely (see Figure \ref{fig:rays_non_transverse}). This contradicts the fact that in such a case we would have replaced the original arc corresponding to this minimizing geodesic with a pair of rays through $p$ instead. Thus there can be no transverse self-intersections in $B_1$.
	\par 
	Suppose we obtain a non-transverse self-intersection. After applying the modified disk flow, $\gamma\cap B_1$ consists of minimizing geodesic arcs connecting either a pair of points of $\partial B_1$ or $p$ and a point of $\partial B_1$. Therefore a non-transverse intersection occurs only when two geodesic segments of $\gamma$ coincide, meet at a point on $\partial B_1$, or meet at $p$. By definition, non-transverse intersections at $p$ are allowed for a $p$-admissible curve. The other two cases occur when two arcs of the original $\gamma\cap B_1$ share either one or both endpoints on $\partial B_1$, which does not occur because $\gamma$ only has self-intersections at $p$. Thus $\gamma$ remains $p$-admissible after applying the modified disk flow.
\end{proof}
The remaining properties of a single curve under the modified disk flow will be proved directly in a manner similar to the proof for the original disk flow (Theorem 1.8 of \cite{hass1994}).	
Next we prove that the modified disk flow does not increase the lengths of curves.
\begin{figure}[ht]
	\centering
	\includegraphics[width=\textwidth]{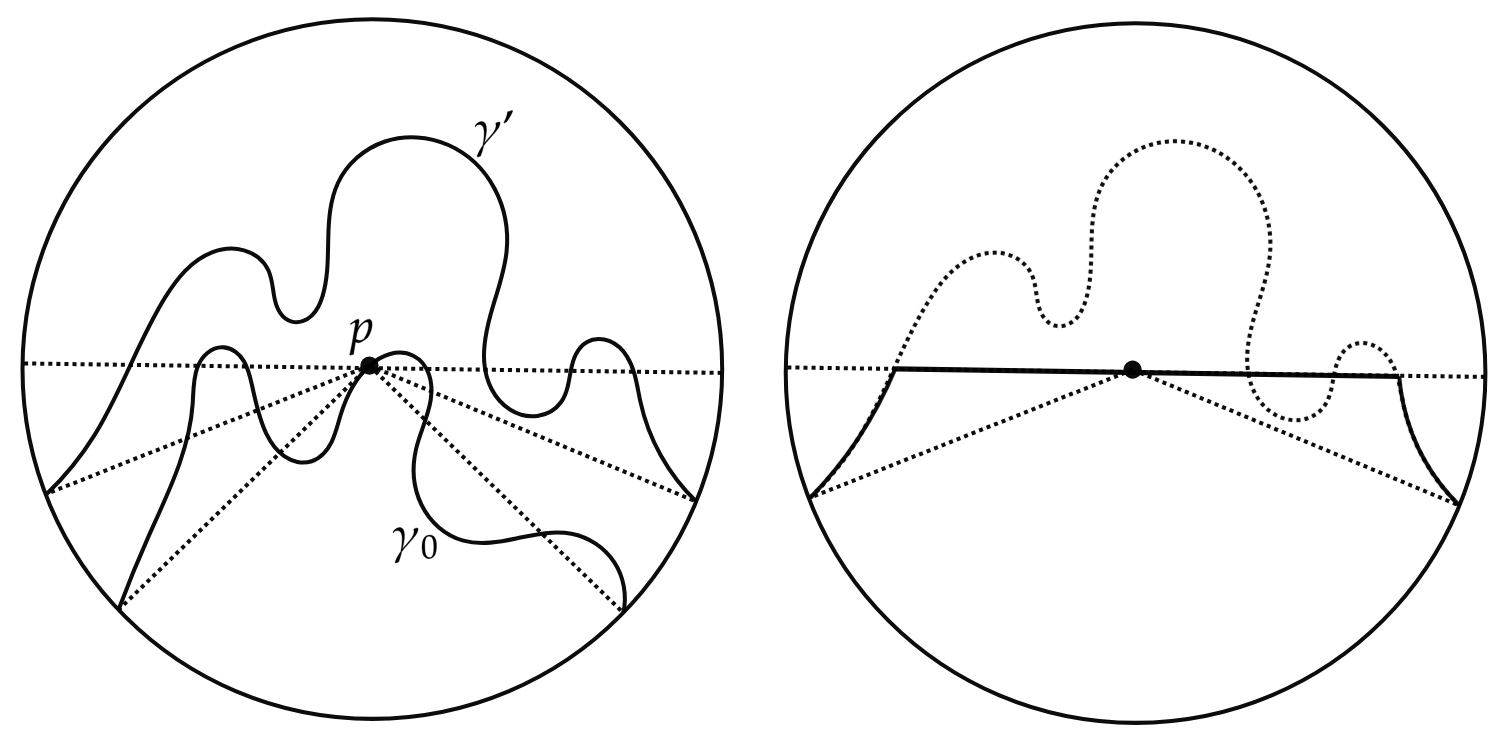}
	\caption{Left: there exists a diameter that lies outside the sector given by the endpoints of $\gamma_0$ and intersects $\gamma'$ at least twice. Right: the curve obtained by replacing a section of $\gamma'$ with a portion of a diameter, which is shorter than $\gamma'$.}
	\label{fig:rays_shorter}
\end{figure}
\begin{lemma}
	\label{lemma:rays_are_shorter}
	Let $\gamma$ be a $p$-admissible curve in $M$. Then the image of $\gamma$ under the modified disk flow is not longer than $\gamma$.
\end{lemma}
\begin{proof} 
	The standard disk flow does not increase lengths since arcs are always being replaced by minimizing geodesics with the same endpoints. Thus the only way we might make $\gamma$ longer is when we replace an arc in $B_1$ with two rays through $p$ that do not form a minimizing geodesic. Any arc of $\gamma\cap B_1$ that passes through $p$ is not lengthened by replacing it with the two rays that connect its endpoints to $p$. Therefore we only need to consider an arc $\gamma'$ with endpoints on $\partial B_1$ that does not pass through $p$. 
	\par
	Suppose that under the modified disk flow $\gamma'$ gets replaced by two rays $\alpha',\beta'$ that do not form a minimizing geodesic. Let $\gamma_0$ be the arc of $\gamma\cap B_1$ containing $\gamma(0)$. As in Lemma \ref{lemma:rays_non_transverse}, $\alpha'$ and $\beta'$ bound a convex sector of angle strictly less than $\pi$ that contains the endpoints of $\gamma_0$. Therefore we can pick a minimizing geodesic $\eta$ that passes through $p$ and connects two points on $\partial B_1$ that lie outside of this sector. Necessarily, the endpoints of $\gamma'$ and $\gamma_0$ lie on the same side of $\eta$. The arcs $\gamma_0$ and $\eta$ both pass through $p$, while $\gamma'$ does not pass through $p$ but also does not cross $\gamma_0$. Therefore $\gamma'$ must intersect  $\eta$ at least twice (see Figure \ref{fig:rays_shorter}, left). Therefore replacing an arc of $\gamma'$ with an arc of $\eta$ that passes through $p$ creates a curve no longer than $\gamma'$ that passes through $p$ (see Figure \ref{fig:rays_shorter}, right). We have now reduced to our preliminary case, and the result follows.
\end{proof}
We now show that a subsequence of the curves obtained by applying the modified disk flow will converge to either a point curve or some number of simple geodesic loops based at $p$. Again, we follow the proof structure for the corresponding theorem for the original disk flow.
\begin{lemma}\label{lemma:modified_disk_flow}
	Let $\gamma$ be a $p$-admissible loop and let $\gamma_t$, $t\in[0,\infty)$, denote the family obtained by applying the modified disk flow. Then the following hold.
	\begin{enumerate}
		\item 
		A subsequence of the curves $\gamma_i$, $i\in\bN$, either converges to the point curve $p$ or to a union of simple geodesic loops $\gamma_\infty$ at $p$. In the latter case, $L(\gamma_\infty)=\lim\limits_{i\to\infty}L(\gamma_i)\leq L(\gamma)$.
		\item 
		$L(\gamma_{i+1})=L(\gamma_{i})$ for some $i\in\mathbb{N}$ if and only if $\gamma_i$ is a point curve or a union of simple geodesic loops at $p$.
        \item 
        Given $i\in\bN$ with $i>0$, there is a number $k_i$, depending only on $\gamma_0$ and the metric ball cover associated with the modified disk flow, such that $\gamma_i$ is a piecewise geodesic with at most $k_i$ breaks.
	\end{enumerate}
\end{lemma}
\begin{proof}
	Lemma \ref{lemma:rays_are_shorter} proves that the lengths of the curves $\gamma_i$ is non-increasing, and hence uniformly bounded by $L(\gamma)$. Therefore we can apply the Arz\`ela--Ascoli theorem to conclude the existence of a convergent subsequence with uniform limit $\gamma_\infty$ satisfying $L(\gamma_\infty)\leq\liminf_{t\to\infty}L(\gamma_t)$. If $\liminf_{t\to\infty}L(\gamma_t)$ is zero, then $\gamma_\infty$ is necessarily a point curve and $L(\gamma_\infty)=\lim\limits_{i\to\infty}L(\gamma_i)$ as claimed. Otherwise, assume $0<\liminf_{t\to\infty}L(\gamma_t)$. We will first show that $\gamma_\infty$ is a piecewise geodesic curve with vertices only at $p$ (i.e., a union of geodesic loops at $p$). 
    \par 
    Suppose this is not the case. Then $\gamma_\infty$ contains a non-geodesic subarc $\sigma$ that is not formed by two minimizing geodesic rays through $p$. 
    We claim that $\gamma_\infty$ can be strictly shortened by applying the modified disk flow. 
    There is some ball $B_j$ that contains a non-geodesic portion of $\sigma$ in its interior. 
    If $\sigma$ is still non-geodesic by the time we apply the shortening procedure in such a ball, then replacing each arc of $\sigma$ in that ball by a minimizing geodesic strictly shortens $\sigma$, and hence strictly shortens $\gamma_\infty$. The other possibility is that $\sigma$ is formed by a pair of geodesics $\sigma_1$ and $\sigma_2$ that meet at a vertex $x$ and the modified disk flow replaces, say, $\sigma_1$ by a geodesic of equal length that meets $\sigma_2$ at angle $\pi$. This replaces $\sigma$ by a geodesic of the same length. Suppose this first occurs when shortening in the ball $B_k$. Then necessarily $x\in B_k$ and $x$ is connected by some arc of $\gamma_\infty\cap B_k$ to a point $y\in \partial B_k$ that lies on the geodesic line covering $\sigma_2$. However, this connecting arc cannot be a geodesic, since it contains a portion of $\sigma_1$ incident to $x$. Therefore when shortening in $B_k$, we are once again replacing a non-geodesic by a minimizing geodesic arc of strictly shorter length.
	\par 
	Thus, applying the modified disk flow to $\gamma_\infty$ produces a strictly shorter curve. 
    We will show that this contradicts the fact that $L(\gamma_\infty)\leq\liminf_{t\to\infty}L(\gamma_t)$, proving that $\gamma_\infty$ is a union of geodesic loops at $p$. By the above argument, there is at least one non-geodesic arc of $\gamma_\infty$ that is not formed by two minimizing geodesic rays through $p$ and that is replaced by a minimizing geodesic under the modified disk flow in some ball $B_j$.   
    Consider the earliest ball $B_i$ where such a replacement occurs. Choose some non-geodesic arc $\sigma_\infty$ of $\gamma_\infty\cap B_i$. By the above argument, no changes to $\gamma_\infty$ can occur before reaching $B_i$. 
    Therefore there is a disk $V\subset B_i$ that is not contained in any earlier ball in the cover and that intersects $\sigma_\infty$ in a non-geodesic arc. 
    For every sufficiently large $k$ (say, $k>K$), the curve $\gamma_k$ has a corresponding arc $\sigma_k$
    that intersects $V$ in a non-geodesic arc. 
    Define $\delta_k>0$ as the difference in length between $\sigma_k\cap V$ and the minimizing geodesic connecting the endpoints of $\sigma_k\cap V$, and let $\delta=\inf_{k>K}\delta_k$. 
    Note that $\delta>0$ for sufficiently large $K$, since the $\sigma_k$ converge uniformly to the non-geodesic arc $\sigma_\infty$. After shortening in the first $i-1$ disks, $\sigma_\infty\cap V$ remains unchanged and each $\sigma_k\cap V$ remains close to $\sigma_\infty\cap V$. Therefore each curve $\gamma_k$ is shortened by at least $\delta$ when applying the disk flow in $B_i$, since the arc containing $\sigma_k\cap V$ is replaced by a minimizing geodesic that is at least $\delta$ shorter in length.
    \par
    Now, by definition we have $L(\gamma_k)-\delta\geq L(\gamma_{k+1})$. Conversely, because lengths are converging we have $ L(\gamma_k)-\delta<L(\gamma_\infty)$ if $k$ is sufficiently large. In combination, we have 
    $
    L(\gamma_{k+1})\leq
    L(\gamma_k)-\delta<L(\gamma_\infty),
    $
    contradicting the fact that $L(\gamma_\infty)\leq\liminf_{t\to\infty}L(\gamma_t)$.
	Thus $\gamma_\infty$ must be a union of geodesic loops at $p$. 
	Note that the loops in the limit are also simple. This is because $\gamma_t$ is $p$-admissible for all $t<\infty$, and hence has only non-transverse self-intersections. Therefore the same is true of $\gamma_\infty$. Because geodesic segments never have non-transverse self-intersections, each loop of $\gamma_\infty$ is simple. Of course, $\gamma_\infty$ may have self-intersections at $p$.
	\par			
	We now address the second claim. If $L(\gamma_{i})=0$, then $\gamma_{i}$ is a point curve and hence $\gamma_j$ is a point curve for all $j>i$. If $L(\gamma_{i})>0$ and $\gamma_i$ is not of the desired form, then, as in the proof of the first claim, an application of the modified disk flow will produce a strictly shorter curve. However, this contradicts the fact that $L(\gamma_{i})=L(\gamma_{i+1})$.
    \par 
    We now prove the third claim. Since $\gamma_0$ is $p$-admissible, there is a number $k_0$ such that $\gamma_0$ is a piecewise geodesic with at most $k_0$ breaks.
    We prove by induction that, for $i\in\bN$ with $i\geq0$, the number of breaks in the curve $\gamma_{i+1}$ is bounded in terms of $i$, the metric ball cover, and $\gamma_0$. Suppose $\gamma_{i}$ has at most $k_{i}$ breaks for some $k_i$. Consider what happens to an arc of $\gamma_{i}\cap B_j$ under the modified disk flow in $B_j$. The first possibility is that this arc is replaced by a minimizing geodesic. If the original arc is already a minimizing geodesic, then the arc, and hence the number of breaks, is unchanged. 
    Otherwise, after replacement the arc has at most two breaks and these breaks must occur at its endpoints on $\partial B_j$. Therefore the number of breaks can only increase if the original arc has only one break and this break occurs in the interior of $B_j$. Therefore the largest possible increase in breaks occurs when every arc in $B_j$ has exactly one interior break. The maximum number of such arcs is equal to the number of breaks. Therefore after all such replacements in $B_j$, we have at most doubled the total number of breaks. 
    \par 
    The other case occurs when an arc of $\gamma_i\cap B_j$ is replaced by a pair of geodesic rays through $p$. Necessarily $j=1$. The resulting arc has one interior break at $p$ and at most two breaks at its endpoints on $\partial B_1$. 
    By Lemma \ref{lemma:rays_are_shorter}, an arc can only be replaced by a pair of rays if this would not increase the total length. Every pair of rays through $p$ has total length $2r_1$, where $r_1$ is the radius of $B_1$. Therefore any arc that is replaced by a pair of rays must have length at least $2r_1$. There are at most $\lfloor L(\gamma_i)/(2r_1)\rfloor$ such arcs. Thus replacements by rays in $B_1$ can increase the total number of breaks by at most $\lfloor3L(\gamma_i)/(2r_1)\rfloor\leq 3L(\gamma_0)/(2r_1)$.
    \par
    The above possibilities are summarized as follows. After shortening in $B_1$, $\gamma_i$ has at most $2k_i+ 3L(\gamma_0)/(2r_1)$ breaks. In each successive ball, we can at most double the number of breaks. After shortening in every ball once, we obtain the curve $\gamma_{i+1}$. Therefore this curve has at most
    $$k_{i+1}= 2^{n-1}(2k_i+ 3L(\gamma_0)/(2r_1))
    $$ 
    breaks, where $n$ is the number of metric balls in our cover.
    By inducting on $i$, we see that $k_{i+1}$ depends only on $\gamma_0$ and our choice of metric ball cover.
\end{proof}
From now on, we will let $\gamma_\infty$ denote the limit of a choice of such a subsequence $\gamma_i$ from the statement of Lemma \ref{lemma:modified_disk_flow}. Our ultimate goal is to find two geometrically distinct geodesic loops based at $p$. Thus we want to prove that if $\gamma_\infty$ is non-trivial, it is not given by iterating a single loop multiple times. This is shown in the following two lemmas.
\begin{lemma}
    \label{lemma:modified_disk_flow_no_doubling}
	Let $\gamma$ be a $p$-admissible loop and let $\eta$ be a simple geodesic loop. If $\gamma_\infty$ is a concatenation of multiple copies of $\eta$ and $-\eta$, then in fact either $\gamma_\infty=\eta^a$ or $\gamma_\infty=-\eta^a$ for some $a\geq 1$.
\end{lemma}
\begin{proof}
	Without loss of generality, orient $\eta$ so that the first loop of $\gamma_\infty$ is $\eta$. We will show that no arc of $\gamma_\infty$ is given by $\eta*-\eta$, and hence that $\gamma_\infty=\eta^a$ for some $a\geq 1$. Suppose instead that some arc converged to $\eta*-\eta$. Consider an arc of the form $\eta\mid_{[0,s]}*-\eta\mid_{[0,s]}$ for some small $s$. No such arc is geodesic and no such arc is arbitrarily well-approximated by a geodesic. Therefore the arcs of $\gamma_i$ that converge to $\eta\mid_{[0,s]}*-\eta\mid_{[0,s]}$ (or at least the portion lying in $B_1\setminus\cup_{i>1}B_i$) must consist of pairs of geodesic rays through $p$. However, the fact that they are rays at all means that they must bound a convex sector that encloses the arc of $\gamma_i\cap B_1$ passing through $\gamma_i(0)$. Let this sector have angle $\theta_i$. Necessarily, $\theta_i$ must tend to 0, which we will show leads to a contradiction. We will consider two cases: when the final loop of $\gamma_\infty$ is $\eta$, and when it is $-\eta$. 
	\par 
	If the final loop of $\gamma_\infty$ is $\eta$, then the angle at $\gamma_\infty(0)$ is greater than zero since it equals the angle of $\eta$ at $p$. Therefore for sufficiently large $i$ the angle at $\gamma_i(0)$ is bounded from below by some $c>0$ (see Figure \ref{fig:doubled}, left). Then $\theta_i$ must be at least as large as $c$ because the corresponding sector contains the arc of $\gamma_\infty$ through $p$. In particular, $\theta_i$ cannot tend to zero, which is a contradiction. 
	\par 
	On the other hand, if the final loop of $\gamma_\infty$ is $-\eta$, then the angle at $\gamma_\infty(0)$ is zero. In this case, consider an arc of $\eta*-\eta$ in $B_1$ of the form $\eta\mid_{[r,1]}*-\eta\mid_{[r,1]}$ for some small $r>0$. As before, this arc is not geodesic and so the arcs of $\gamma_i$ that converge to it must also consist of pairs of geodesic rays through $p$ whose angles at $p$ tend to 0 (see Figure \ref{fig:doubled}, right). Therefore these rays must also bound an acute sector containing the arc of $\gamma_\infty$ through $\gamma_\infty(0)$. However, this sector is disjoint from the previous one by $p$-admissibility, so they cannot both contain the arc of $\gamma_\infty$ through $\gamma_\infty(0)$. Thus we once again obtain a contradiction.
\end{proof}	
\begin{figure}[ht]
	\centering
	\includegraphics[width=\textwidth]{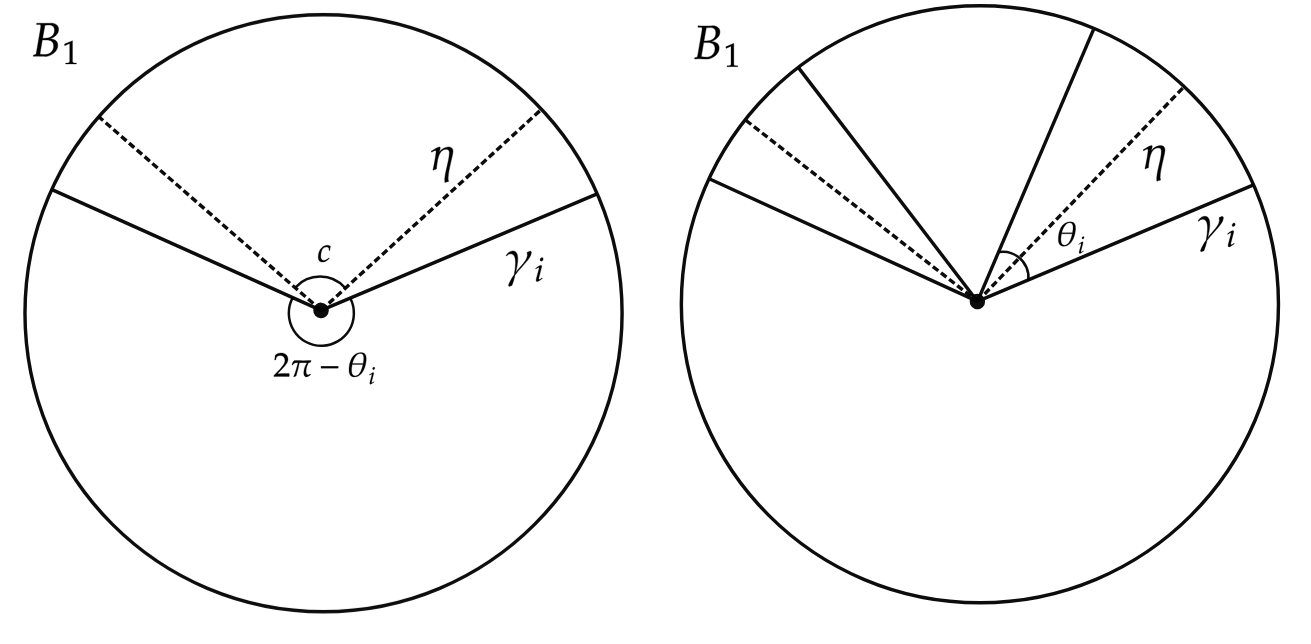}
	\caption{Left: the sector bounded by $\gamma_i$ when the final loop of $\gamma_\infty$ is $\eta$. Right: the two sector bounded by $\gamma_i$ when the final loop of $\gamma_\infty$ is $-\eta$.}
	\label{fig:doubled}
\end{figure}
\begin{lemma}
	\label{lemma:modified_disk_flow_prime}
	Let $\gamma$ be a $p$-admissible loop. If $\gamma_\infty=\eta^a$ for some geodesic loop $\eta$, then $a=1$ (i.e., $\gamma_\infty$ is prime).
\end{lemma}
\begin{proof}
	Let $r$ be small enough that $B_{r}(p)$ is contained in $B_1\setminus \cup_{i>1}B_{i}$ and $B_r(p)\cap \eta$ consists of exactly one arc. This is possible because $\eta$ is simple by Lemma \ref{lemma:modified_disk_flow}. Pick $\epsilon<r/2$ and let $i$ be large enough that $\gamma_i$ is within distance $\epsilon/2$ of $\gamma_\infty$ in the sup norm. 
	Because $\gamma_i$ is $p$-admissible, we can apply a homotopy of very small width to make $\gamma_i$ into a simple curve that still lies in an $\epsilon$-neighbourhood of $\eta$. Consider the arcs of $\gamma_i\cap B_{r}(p)$. We will call these arcs $\Gamma_i^j$, $j\in\{1,\cdots,k\}$, ordered as follows. Letting $\partial B_{r}(p)$ be positively oriented, choose some arc of $\gamma_i\cap B_{r}(p)$ that, with $\partial B_{r}(p)$, cobounds a region of $B_{r}(p)$ whose interior contains no other arc of $\gamma_i$. This is possible due to the fact that $\gamma_i$ is simple. Call this arc $\Gamma_i^1$. Number the other arcs of $\gamma_i\cap B_{r}(p)$ in the order in which they are encountered when following $\partial B_1$ from the endpoint of $\Gamma_1^i$. This ordering is well-defined because $\gamma_i$ has no transverse self-intersections and self-intersects only at $p$. Moreover, let $x_j$ be the initial point of the arc $\Gamma_i^j$ and let $y_j$ be the end point, according to the orientation of $\gamma_i$. These points are necessarily arranged along $\partial B_1$ in the cyclic order $x_1, x_2, x_3,\cdots, x_k,y_k,y_{k-1},\cdots, y_2, y_1$.
	\par 
	\begin{figure}[ht]
		\centering
		\includegraphics[width=\textwidth]{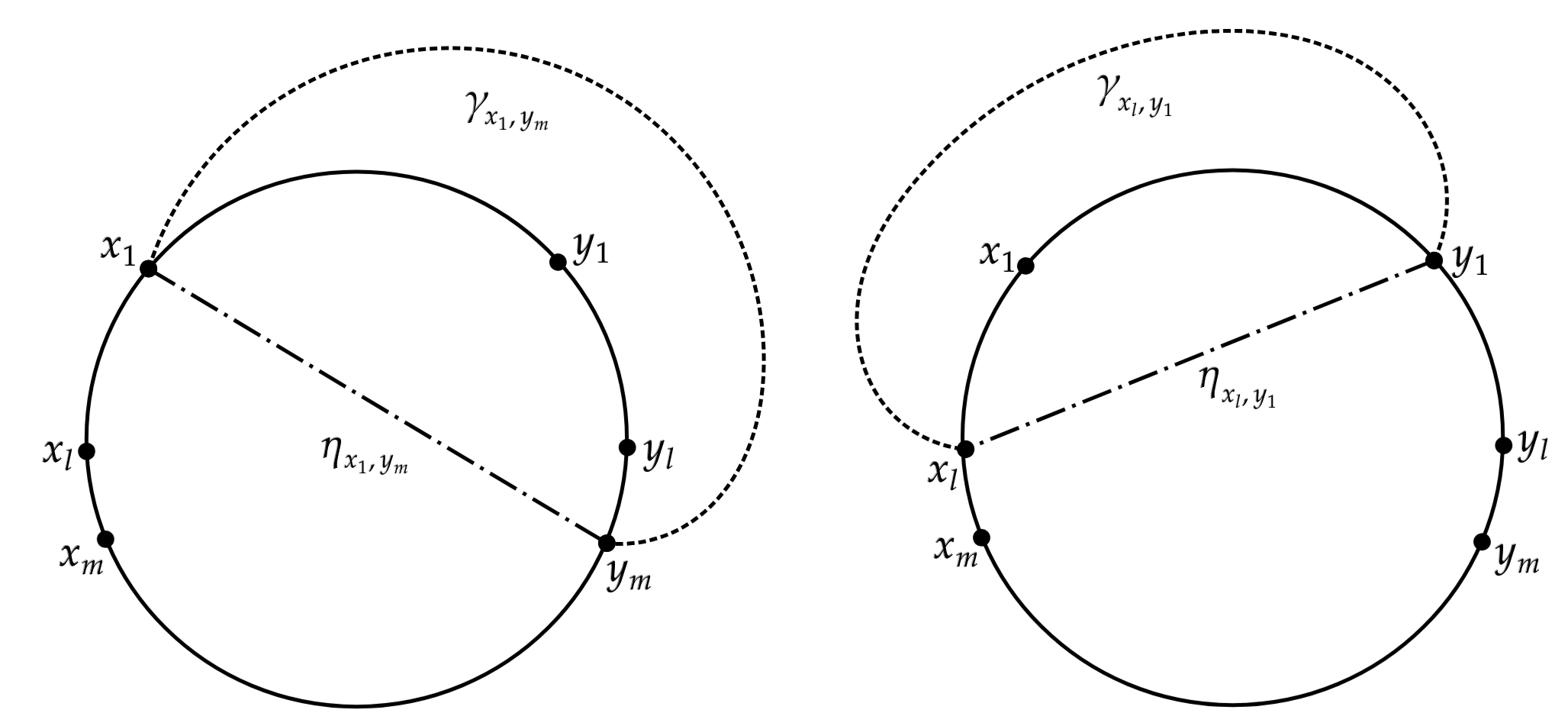}
		\caption{Separating the endpoints of arcs $\Gamma_i^j$ by $\gamma_{x_1,y_m}\cup\eta_{x_1,y_m}$ (left) and  $\gamma_{x_l,y_1}\cup\eta_{x_l,y_1}$ (right).}
		\label{fig:crossing}
	\end{figure}	
	Because $B_{r}(p)\cap \eta$ only contains the arc of $\eta$ passing through $p$, for sufficiently small $\epsilon$ every $\Gamma_i^j$ represents an arc of $\gamma_i$ that is converging to an arc of the form $\eta\mid_{[1-s,1]}*\eta\mid_{[0,s]}$ for some $s$. Moreover, every arc of $\gamma_i$ that converges to a copy of $\eta$ corresponds to exactly one arc $\Gamma_i^j$.
	We will now prove that $k=1$, and hence that $\gamma_\infty = \eta$. 
	Following $\gamma_i$ forward from $y_1$, we next encounter some $x_l$. Let the arc of $\gamma_i$ joining these two points be denoted $\gamma_{y_1,x_l}$, noting that this arc lies entirely outside of $B_1$. Going backwards along $\gamma_i$, $x_1$ is followed by some $y_m$. Let the arc of $\gamma_i$ joining these two points be denoted $\gamma_{y_m,x_1}$, which we note again lies outside $B_1$. If $\rho_{x_l,y_1}$ is the minimizing geodesic in $B_1$ connecting $x_l$ and $y_1$, then $x_1$ is separated from any $y_i$ with $i>1$ and any $x_j$ with $j>l$ by the closed curve $\gamma_{y_1,x_l}*\rho_{x_l,y_1}$ (see Figure \ref{fig:crossing}, left). Similarly, $y_1$ is separated from any $x_i$ with $i>1$ and any $y_j$ with $j>m$ by the closed curve $\gamma_{y_m,x_1}*\rho_{x_1,y_m}$ (see Figure \ref{fig:crossing}, right). However, this means that the arcs $\gamma_{y_1,x_l}$ and $\gamma_{y_m,x_1}$ must either cross, which they cannot by $p$-admissibility, or they coincide (i.e., they are the same portion of $\gamma$). Thus we must have $x_l=x_1$ and $y_m=y_1$, meaning that $k=1$ and hence that $\gamma_\infty=\eta$.
\end{proof}	
\noindent 
Consequently, we can find a subsequence of curves with decreasing lengths that converges to either a point or a union of simple geodesic loops. If we obtain loops, then either at least two of our loops have distinct images or we obtain a single prime loop. Note that we cannot exclude the possibility of obtaining a union of multiple distinct geodesic loops, as such a set of curves is fixed by the modified disk flow for a suitable parameterization.

\subsection{Short Meridional Slicings}
\label{sec:short_meridional_slicings}

We now need to find suitable families of curves to which we can apply the modified disk flow. In the standard Lusternik--Schnirelmann proof, one can use cycles made of meridians on the round sphere. While we could simply take the images of those curves under a diffeomorphism, this would not provide us with a length bound. Instead, we do the following. If $p$ is the point where we want to base our geodesic loops, pick a point $q$ that is at locally maximal distance from $p$ in $M$. These two points are analogous to the poles on the sphere. By Berger's Lemma, for any $v\in T_qM$ there is a minimizing  geodesic $\gamma$ from $q$ to $p$ with $\langle \gamma'(0),v\rangle \geq 0$. Therefore we can pick either two or three such geodesics, denoted $\tau_i$ and cyclically indexed, such that $\langle \tau_i'(0),\tau_{i+1}'(0) \rangle \geq 0$. Our meridional slicing will be given by interpolating between these geodesics through simple curves that do not mutually internally intersect. We call such a family a strictly monotone meridional slicing, which we define below. We first define strict monotonicity, which we will also apply to families of curves that are not meridional slicings.
\begin{definition}[Strict Monotonicity]
    \label{def:monotonicity}
A free homotopy $t\mapsto f_t$ is called strictly monotone if every $f_t$ is simple and no pair of curves $f_t$ and $f_s$ with $s\not=t$ intersect. A fixed-endpoint homotopy $t\mapsto f_t$ is called strictly monotone if every $f_t$ is simple and no pair of curves $f_t$ and $f_s$ with $s\not=t$ intersect except at their endpoints.
\end{definition} 
\begin{definition}[Meridional Slicings]
    \label{def:slicing}
    A map $t\mapsto\Gamma_t$ from $S^1$ into the space of curves on $M$ is called a meridional slicing if the following properties hold.
    \begin{enumerate}
        \item 
        There exist points $p$ and $q$ in $M$ such that $\Gamma_t(0)=p$ and $\Gamma_t(1)=q$ for all $t\in S^1$.
        \item 
        $\Gamma$ induces a map on $M$ of non-zero degree.		
    \end{enumerate}
    We call $\Gamma_t$ strictly monotone if it is strictly monotone when regarded as a fixed-endpoint homotopy-- that is, if every $\Gamma_t$ is simple and $\Gamma_t\cap\Gamma_s=\{p,q\}$ for all distinct $t,s\in S^1$.
\end{definition} 
\par 
We will construct these interpolations in two steps. First, we will contract the curves $-\tau_i*\tau_{i+1}$ through loops based at $p$ whose lengths are bounded in terms of the diameter $d$. We then apply the following lemma to convert this homotopy to a homotopy between $\tau_i$ and $\tau_{i+1}$ through curves whose lengths are also bounded. The argument is as follows.
If the homotopy from $\gamma_1*-\gamma_2$ to $p$  given in the statement of the lemma is denoted by $\Gamma_t$, then the homotopy from $\gamma_1$ to $\gamma_2$ is given explicitly by applying the homotopy $-\Gamma_{t}*\gamma_1$, $t\in[0,1]$, followed by the homotopy $\gamma_2*-\gamma_1\mid_{[0,1-t]}*\gamma_1\mid_{[0,1-t]}$, $t\in[0,1]$.
\begin{lemma}
	\label{lemma:maeda}
	Suppose $\gamma_1$ and $\gamma_2$ are arcs with shared endpoints $p$ and $q$. Suppose also that $\gamma_1*-\gamma_2$ can be homotoped to $p$ via loops based at $p$ of length at most $L$. Then $\gamma_1$ can be homotoped to $\gamma_2$ through curves of length at most $L(\gamma_1)+L$.\qed
\end{lemma}
First, however, we need to determine how to contract the loops $-\tau_i*\tau_{i+1}$. If we use the modified disk flow, fixing $p$, the homotopy may stall at a union of geodesic loops instead of a point curve. Indeed, it is possible there are only two geodesics, $\tau_1,\tau_2$, which cover two halves of a single closed geodesic through $p$. Such a curve will be fixed by the (modified) disk flow. Therefore a more involved algorithm is needed.
\par 
For brevity, we utilize the following definition.
\begin{definition}[Geodesic Digons]
	A closed region bounded by a closed curve consisting of the concatenation of two (minimizing) geodesic segments is called a (minimizing) geodesic digon. The two endpoints of the segments are called vertices.
\end{definition}
\noindent
In particular, the loops we want to contract are the boundaries of minimizing geodesic digons with vertices at $p$ and $q$, and with an acute angle at $q$. The behaviour of the angles at $p$ and $q$ will help us control the behaviour of the modified disk flow.
\par
\begin{figure}
	\centering
	\includegraphics[width=\textwidth]{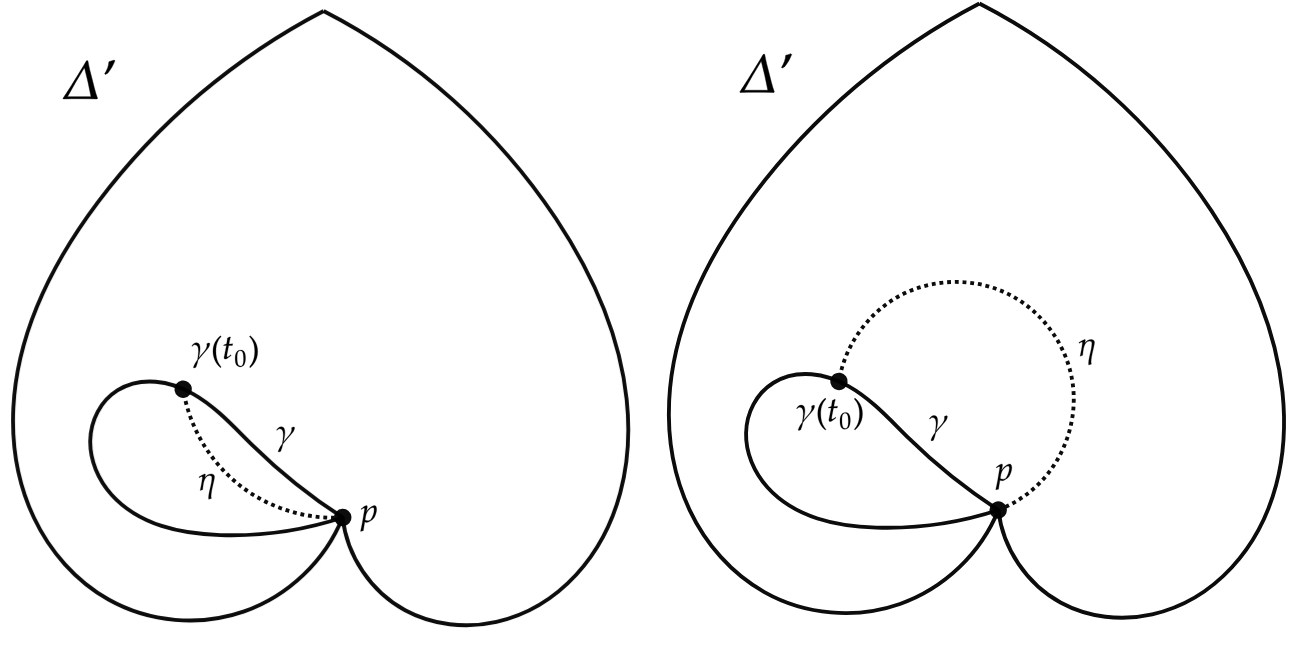}
	\caption{The two possible ways a minimizing geodesic $\eta$ (dashed) can intersect the geodesic loop $\gamma$ within the digon $\Delta'$.}
	\label{fig:gamma_sep}
\end{figure}

 
We will also use the fact that $q$ lies in the cut locus of $p$, allowing us to exploit the structure of the cut locus. In particular, we can take advantage of the following theorem.	
\begin{theorem}[S. B. Myers, 1935 \cite{myers1935}]
	\label{theorem:myers}
	Suppose $M$ is a Riemannian 2-sphere with an analytic metric. Then for any $x\in M$, the cut locus of $x$ consists of either a single point or a finite embedded tree. If the cut locus is a single point, this point is conjugate on every geodesic connecting it to $x$. If instead the cut locus is a tree, then the degree of any point $y$ equals the number of minimizing geodesics connecting $x$ to $y$. Moreover, $M$ consists of a disk with geodesic polar co-ordinates centred at $x$ glued along the cut locus of $x$.
	\qed
\end{theorem}	
\noindent 
Because the cut locus is finite, we will be able to induct on the number of vertices, say $v$, that are contained in the digon bounded by $-\tau_{i}*\tau_{i+1}$, which we will call $\Delta$. Using Theorem \ref{theorem:myers}, we will show that either $q$ is at a vertex of the cut locus of $p$, or we can homotope $\Delta$ to a digon that contains at most $v-1$ vertices of the cut locus. In the latter case, we proceed by induction on $v$. If $q$ is at a vertex, then some number of edges of the cut locus incident to $q$ lie in $\Delta$ (or else $\Delta$ would be easily contractible to $p$). The number of these branches corresponds to the number of minimizing geodesic digons with vertices at $p$ and $q$ that $\Delta$ covers. We then consider a digon lying in $\Delta$, say $\Delta'$, with the property that the angle of $\Delta'$ at $q$ is at most $\pi$ and one edge of $\Delta'$ coincides with one edge of $\Delta$. Since they share an edge, contracting $\Delta'$ to a point will allow us to homotope $\Delta$ to the smaller digon $\Delta\setminus \Delta'$. We therefore apply the modified disk flow to the boundary of $\Delta'$. By the following lemma, any loops obtained from applying the modified disk flow lie within $\Delta'$.
\begin{lemma}
	\label{lemma:shortening_digons}
	Suppose $\Omega$ is a closed region bounded by a piecewise geodesic curve $\partial\Omega$ that passes through $p$. Suppose every interior angle at every vertex of $\partial\Omega$, except perhaps the one at $p$, is at most $\pi$. Then any simple curve $\gamma$ in $\Omega$ (including $\partial\Omega$ if it is simple) remains in $\Omega$ under the modified disk flow. Consequently, the resulting family consists of simple curves that do not mutually intersect transversely.
\end{lemma}
\begin{proof}
	Let the family obtained by applying the modified disk flow to $\gamma$ be denoted by $\gamma_t$. We need to ensure that for every $t$, $\gamma_t$ lies within $\Omega$. It is enough to check that $\gamma$ remains in $\Omega$ after applying the disk flow once. In each ball $B_i$ with $i>1$, $\Omega\cap B_i$ is a convex region containing $\gamma\cap B_i$. Therefore when each arc of $\gamma\cap B_i$ is replaced by a minimizing geodesic segment, this segment will necessarily lie in $\Omega\cap B_i$. 
	\par 
	Now consider $B_1$. If the interior angle of $\partial\Omega$ at $p$ is at most $\pi$, then $\Omega\cap B_1$ is also convex and the result holds. 
	Suppose instead the angle at $p$ is greater than $\pi$. As before, let $\gamma_0$ denote the arc of $\gamma\cap B_1$ through $\gamma(0)$ and denote the two rays connecting $p$ to the endpoints of $\gamma_0$ by $\alpha$ and $\beta$. Consider the closed sector $S$ in $\Omega\cap B_1$ bounded by $\alpha$ and $\beta$ (since $p\in \partial\Omega$, there is a unique such sector).
	Then $(\Omega\cap B_1)\setminus S$ consists of two other sectors $S'$ and $S''$ with vertex at $p$.
	Suppose there is another arc of $\gamma\cap B_1$, say $\gamma'$.  Because $\gamma$ is simple, $\gamma'$ does not pass through $p$ and hence must lie in exactly one of the three sectors.
	As in Lemma \ref{lemma:rays_acute_angle}, $\gamma'$ bounds a fourth sector with vertex at $p$ that does not contain the endpoints of $\gamma_0$, and $\gamma'$ is only replaced by a minimizing geodesic if that sector has an acute angle at $p$ (i.e., if it is convex). If $\gamma'$ is replaced by two geodesic rays through $p$, then its image under the modified disk flow will necessarily lie in $\Omega$. Otherwise, the fourth sector has an acute angle. Because $\gamma'$ does not pass through $p$, this sector is strictly contained in either $S,S'$ or $S''$. Therefore $\gamma'$ lies in a convex sector within $\Omega$, and hence the minimizing geodesic connecting the endpoints of $\gamma'$ lies entirely within $\Omega\cap B_1$, as desired. 
\end{proof}
\noindent 
Therefore we either succeed in contracting $\partial \Delta'$ or we encounter some simple geodesic loops in $\Delta'$ based at $p$ of length at most $L(\partial \Delta')\leq 2d$. If $M$ contains two such loops, then Theorem \ref{theorem:main} is proved. Therefore we can assume that either we can contract $\partial \Delta'$ or $\Delta'$ contains the only such loop. 
Then any digon satisfying the hypotheses of Lemma \ref{lemma:shortening_digons} that lies outside of $\Delta'$ is contractible via the disk flow. Therefore if we can contract the loop in $\Delta'$, we can homotope $\Delta$  to $\Delta\setminus\Delta'$ via Lemma \ref{lemma:maeda}. Then $\Delta\setminus\Delta'$, and hence $\Delta$ itself, is easily contractible. We detail this process in the following lemma. 
\begin{lemma}
	\label{lemma:digons_and_loops} Suppose the metric on $M$ is analytic.
	Let $\Delta$ be a minimizing geodesic digon with vertices $p$ and $q$. Suppose also the interior angle of $\partial \Delta$ at $q$ is at most $\pi$. Then either there are at least two simple geodesic loops at $p$ of length at most $2d$ in the interior of $\Delta$ or $\partial \Delta$ can be contracted to $p$ through loops based at $p$ of length at most
	$6d$.
\end{lemma}
\begin{proof}
	Assume that $\Delta$ contains at most one such loop. Recall by Theorem \ref{theorem:myers} that the cut locus of $p$ has the structure of an embedded tree or a single point. We will consider a point to be a tree consisting of a single vertex. Therefore, suppose that the closure of $\Delta$ contains exactly $v$ vertices of the cut locus of $p$. We will proceed by induction on $v$ and the number of loops contained in $\Delta$, which we denote by $m$, to show that $\Delta$ can be contracted through loops based at $p$ of length bounded by some function $c(v,m)$. Note that $v\geq 1$, since $q$ lies on the cut locus of $p$, and $0\leq m\leq 1$ by assumption. 
	\par 
	Note that if $q$ lies within an edge of the cut locus (i.e., it is not a vertex), then $\partial\Delta$ can be homotoped to a minimizing geodesic digon with vertices at $p$ and some cut locus vertex $q'\in\Delta$. This is possible since minimizing geodesics vary continuously as their endpoint moves within an edge of the cut locus. Thus we can assume that $q$ lies on a vertex.
	\par
	Our base case is when $v=1$. Since $q$ is a vertex, the cut locus intersects $\partial \Delta$ only at $q$. Therefore $\partial \Delta$ can be contracted through curves of length at most $c(1,m)=2d$ by perturbing it off of $q$ and applying the modified disk flow.
	\par 
	Suppose $v>1$. Note that the cut locus must be a tree in this case.
	Therefore by Theorem \ref{theorem:myers}, the total number of minimizing geodesic segments connecting $p$ and $q$ depends on the degree of $q$. By the same result, at least one edge of the cut locus incident to $q$ is contained in $\Delta$. If only one edge lies in $\Delta$, then $\Delta$ can be homotoped to a new digon (contained in $\Delta$) by moving its two bounding geodesics along this edge until they reach the next vertex contained in $\Delta$. This new digon contains $v-1$ vertices, so by induction it can be contracted through loops of length $c(v-1,m)$. 
	\par 
	Suppose instead that there are at least two edges of the cut locus incident to $q$ within $\Delta$. Then the geodesics emanating from between these edges split $\Delta$ into minimizing geodesic digons whose interiors are disjoint. Exactly two of these are bounded by one of the edges of $\partial \Delta$ and some other minimizing geodesic. Because the angles between adjacent minimizing geodesics at each vertex must sum to $2\pi$, at least one of these two digons has an angle at $q$ that is at most $\pi$. Call this digon $\Delta'$. 
	Apply the modified disk flow to $\partial\Delta'$. 
	\par
	If $\partial\Delta'$ contracts to a point, then we can apply Lemma \ref{lemma:maeda} to homotope one of the edges of $\partial \Delta'$ to the other via curves of length at most $3d$. Therefore $\partial \Delta$ can be homotoped to $\partial (\Delta\setminus\Delta')$ though curves of length at most $4d$, since $\partial \Delta'$ and $\partial\Delta$ share an edge. Because $\Delta'$ contains an edge of the cut locus and hence at least one vertex in its interior, $\Delta\setminus\Delta'$ contains at most $v-1$ vertices of the cut locus. Thus by induction this digon can be contracted through loops of length $c(v-1,m)$. Thus $\Delta$ is contractible through loops of length at most $\max\{4d, c(v-1,m)\}$.
	\par 
	If $m=1$, then the other possibility is that the modified disk flow will contract $\partial\Delta'$ to the unique simple geodesic loop based at $p$ of length at most $2d$, which we will call $\gamma$. We need to contract this loop in order to finish our contraction of $\Delta$. Because $\gamma$ is a geodesic loop (and hence minimizing for at most half its length), it must also intersect the cut locus of $p$. Follow $\gamma$ from $\gamma(0)$ to its first intersection point with the cut locus, say $\gamma(t_0)$. By definition, the arc $\gamma\mid_{[0,t_0]}$ is minimizing. 
	If $\gamma\mid_{[t_0,1]}$ is also minimizing, then $\gamma$ itself bounds a minimizing geodesic digon in $\Delta'$ with vertices at $p$ and $\gamma(t_0)$. Moreover, this digon has angle $\pi$ at $\gamma(t_0)$
	and does not contain any short geodesic loops in its interior. Therefore by induction we can contract $\gamma$ through curves of length at most $c(v-1,0)$. Consequently, by Lemma \ref{lemma:maeda} $\partial \Delta$ can be homotoped to $\partial (\Delta\setminus\Delta')$ though curves of length at most $2d+c(v-1,0)$, and then $\Delta\setminus\Delta'$ can be contracted through curves of length at most $c(v-1,0)$.
	\par 
	On the other hand, suppose $\gamma\mid_{[t_0,1]}$ is not minimizing. If $\gamma(t_0)$ is a leaf vertex, then as before we can contract $\gamma$ through curves of length at most $2d$. Otherwise, there is some other minimizing geodesic $\eta$ connecting $p$ to $\gamma(t_0)$. We claim that $\eta$ splits the region enclosed by $\gamma$ into two digons (see Figure \ref{fig:gamma_sep}, left). Because $\eta$ is a minimizing geodesic, the only other potential possibility is that $\eta$ lies outside the region enclosed by $\gamma$ but inside $\Delta'$ (see Figure \ref{fig:gamma_sep}, right). In this case, consider the region bounded by $\eta$, $\gamma\mid_{[t_0,1]}$ and $\partial\Delta$. This region is bounded by a piecewise geodesic with all interior angles acute except perhaps for one of the angles at $p$. Therefore by Lemma \ref{lemma:shortening_digons} the image of $\partial\Delta'$ under the modified disk flow must have remained in this region. This contradicts the fact that the image converged to $\gamma$. Thus $\eta$ must split the region bounded by $\gamma$ into two digons. 
	The angles formed by this pair of digons at $\gamma(t_0)$ are both acute since they lie within the region bounded by $\gamma$. Moreover, neither of these digons contain a simple geodesic loop at $p$ of length at most $2d$. Therefore by Lemma \ref{lemma:shortening_digons} the boundaries of each digon can be shortened to a point through loops of length at most $2d$. As before, we can then contract $\gamma$ (and hence $\Delta'$) through loops of length at most $4d$ and then homotope $\partial \Delta$ to $\partial (\Delta\setminus\Delta')$ though loops of length at most $6d$. Then $\Delta\setminus\Delta'$ can be contracted via curves of length at most $c(v-1,0)$.
	\par 
	In summary, we have $c(1,0)= 2d$ and $c(v,0)= \max\{4d, c(v-1,0)\}$, and hence $c(v,0)= 4d$. Moreover, $c(1,1)= 2d$ and 
	\begin{align*}
		c(v,1)= \max\{6d, c(v-1,1), 2d+c(v-1,0)\}=\max\{c(v-1,1), 6d\}.
	\end{align*} Thus we obtain $c(v,m)= 6d$.
\end{proof}
Now we know how to contract digons, and hence we know how to contract the disks we obtain when applying Berger's lemma. Therefore we are almost ready to construct our strictly monotone meridional slicing. Note that monotonicity is required because we will be forming our curve families from pairs of meridional curves, and we need the resultant curves to be $p$-admissible. The existence of a strictly monotone meridional slicing with endpoints at $p$ and $q$ is given by the following lemma.
\begin{lemma}
	\label{lemma:monotone}
	The contracting homotopy described in Lemma \ref{lemma:digons_and_loops} can be made strictly monotone.
\end{lemma}
\begin{proof}
	The homotopy in Lemma \ref{lemma:digons_and_loops} is given by a concatenation of homotopies of digons and geodesic loops. Because any two minimizing geodesics emanating from $p$ do not intersect except at $p$ and perhaps a point on the cut locus, the homotopy formed by moving a pair of minimizing geodesics down an edge of the cut locus of $p$ is strictly monotone. Moreover, the image of the homotopy lies between the digons formed at the vertices of the edge. Therefore it is enough to check that when shortening a digon in our process, the resulting homotopy is strictly monotone and its image lies within the digon. By Lemma \ref{lemma:shortening_digons}, the homotopy $\gamma_t$ consists of simple curves that do not intersect transversely. The claim is proven by slightly perturbing this homotopy to ensure that the curves do not also mutually intersect non-transversely. 
\end{proof}
\noindent Now we combine the results of this section to form our strictly monotone meridional slicing. 
\begin{lemma}
	\label{lemma:really_short_slicing}
	Given $p \in M$, pick any $q$ at locally maximum distance from $p$ in $M$. Suppose there is at most one geodesic loop based at $p$ of length at most $2d$ in $M$. Then for any $\epsilon>0$ there is a strictly monotone meridional slicing $\Gamma_t$ of $M$ through curves with endpoints at $p$ and $q$ of length at most $7d+\epsilon$.
\end{lemma}
\begin{proof}
	First we find two or three minimizing digons with vertices at $p$ and $q$ as detailed in the introduction to this section. Because $q$ is at locally maximal distance from $p$, Berger's lemma ensures that given any $v\in T_qM$ there is a minimizing geodesic $\gamma$ connecting $q$ to $p$ with $\langle \gamma'(0),v\rangle \geq 0$. In particular, there are either two or three minimizing geodesics $\tau_i$ from $q$ to $p$ of length at most $d$, ordered cyclically, so that each pair $\tau_i*-\tau_{i+1}$ bounds a disk $D_i$ and the closure of these disks cover $M$. If there are only two such disks, they both form angle $\pi$ at $q$. Therefore at least one is convex. If there are three disks, then because the angle formed at $q$ between each adjacent pair is at most $\pi$, at least two of these disks must be convex. If there is a third disk it may have an obtuse angle at $p$. 
	\par 
	Contract each $\partial D_i= -\tau_i*\tau_{i+1}$ through loops based at $p$ via Lemma \ref{lemma:digons_and_loops}, obtaining a homotopy $H^i_t$ with $L(H^i_t)\leq 6d$. By Lemma \ref{lemma:monotone}, we can assume that this homotopy is strictly monotone. Moreover, the curves $H^i_t$ are contained in $D_i$ by Lemma \ref{lemma:shortening_digons}. We now use these homotopies to construct a meridional slicing of $M$, although we will have to be careful to maintain monotonicity. One way to homotope $\tau_i$ to $\tau_{i+1}$ is to follow the curves $\tau_i* H_t^{i}$ from $\tau_i=\tau_i* \{p\}$ to $\tau_i* -\tau_i*\tau_{i+1}$, and then contract $\tau_i* -\tau_i$ along itself as in Lemma \ref{lemma:maeda}. However, we need to alter this homotopy to be strictly monotone.
	\par 
	Pick $\epsilon>0$ small enough that $B_{\epsilon}(p)$ is a totally normal metric ball. Let $f^i(t)$ be the first time that $H_t^i$ intersects $\partial B_{\epsilon (t+1)/2}(p)$, so that $P_t^i=H_t^i\mid_{[f^i(t),1]}$ is a curve from some point on $\partial B_{\epsilon (t+1)/2}(p)$ to $p$. The curve $P^i_t$ will form the second half of the meridian $\Gamma_t$. The first half will be given by a curve $\sigma_t^i$ from $q$ to $H_t^i(f^i(t))$, defined so that the curves $\sigma_t^i$ do not mutually intersect except at $q$ and do not intersect any $H_s^i$ (except perhaps at their endpoints if $s=t$). 
	\par 
	First we define a curve $T_t^i$ from $q$ to $\partial B_{\epsilon}(p)$, which will form the first portion of $\sigma_t^i$. Because $\tau_i$ is minimizing between $p$ and $q$, $\tau_i\cap \partial B_{\epsilon}(p)$ does not lie on the cut locus of $q$. Therefore, since the cut locus of $q$ is a tree (or a mere point), there is some arc $\Sigma^i$ of $\partial B_{\epsilon}(p)\cap D_i$ starting at $\tau_i\cap \partial B_{\epsilon}(p)$ that does not lie in the cut locus of $q$. Parameterize this arc by $[0,1]$ and let $T_t^i$ be the minimizing geodesic connecting $q$ to $\Sigma^i(t)$, ordered so that $T_0^i$ is an arc of $\tau_i$ (see Figure \ref{fig:example_hom}, left). Necessarily, the curves $T_t^i$ are simple, lie in $D_i\setminus B_{\epsilon}(p)$ and do not mutually intersect except at $q$. 
	\par
	Now we connect $T_t^i$ to $H_t^i(f^i(t))$. Let $\rho^i_t:[0,1]\to M$ be a portion of a radial geodesic in $B_{\epsilon}(p)$ that connects $T_t^i(1)\in \partial B_{\epsilon}(p)$ to a point on $\partial B_{\epsilon (t+1)/2}(p)$ (see Figure \ref{fig:example_hom}, centre). Let $\eta^i_t:[0,1]\to M$ be the minimizing geodesic that starts at $\rho^i_t(1)$ and ends at the first point of $\partial B_{\epsilon (t+1)/2}(p)\cap H_t^i$ it meets (see Figure \ref{fig:example_hom}, right). The curve $\sigma_t^i$ is then defined as $\sigma^i_t=T_t^i*\rho_t^i*\eta_t^i$. We stress that, again, these curves are simple and do not mutually intersect except at their endpoints.
	\par
	We are now ready to build the portion of $\Gamma_t$ given by homotoping $\tau_i$ to $\tau_{i+1}$. Recall that $P^i_t=H_t^i\mid_{[f^i(t),1]}$ and fix a small $\delta>0$. First, homotope $\tau_i$ to $\sigma^i_\delta*P^i_{\delta}$ via some strictly monotone homotopy, choosing $\delta$ and $\epsilon$ so that this homotopy passes through reasonably short curves. Then use the homotopy $t\mapsto \sigma^i_t*P^i_t$ to pass from $\sigma^i_\delta*P^i_{\delta}$ to $\sigma^i_{1-\delta}*P^i_{1-\delta}$. For small enough $\delta$ and $\epsilon$, the curve $\sigma^i_{1-\delta}*P^i_{1-\delta}$ is very close to $\tau_i*-\tau_i*\tau_{i+1}$. Therefore we can apply some strictly monotone homotopy through reasonably short curves to contract one ``arm'' of $\sigma^i_{1-\delta}*P^i_{1-\delta}$ (the portion approximating $\tau_i*-\tau_i)$, leaving a curve that lies very close to $\tau_{i+1}$. We then homotope this curve to $\tau_{i+1}$.
	Our desired strictly monotone meridional slicing $\Gamma_t$ is given by using these homotopies to pass from $\tau_1$ to $\tau_2$, then to $\tau_3$ if it exists, and then back to $\tau_1$. Given any constant $c>0$, we can choose $\epsilon$ and $\delta$ small enough that the curves $\Gamma_t$ are bounded in length by $d+6d+c$. 
    Moreover, because the homotopies between $\tau_i$ and $\tau_{i+1}$ are strictly monotone and lie in disjoint disks $D_i$, this family represents a degree $\pm1$ map into $M$, and hence is indeed a meridional slicing.
\end{proof}	
We obtain the curves $H_t^i$ by applying the disk flow a finite number of times to curves that are either geodesics or a union of two geodesics. Therefore Property 3 of Lemma \ref{lemma:modified_disk_flow}, which says that each application of the (modified) disk flow can only increase the number of breaks by a finite amount, proves that each $H_t^i$ is a piecewise geodesic with at most $k$ breaks for some $k>0$. Therefore the curves in the meridional slicing obtained above are also piecewise geodesics with a universally bounded number of breaks.
\begin{figure}[ht]
	\centering
    \includegraphics[width=\textwidth]{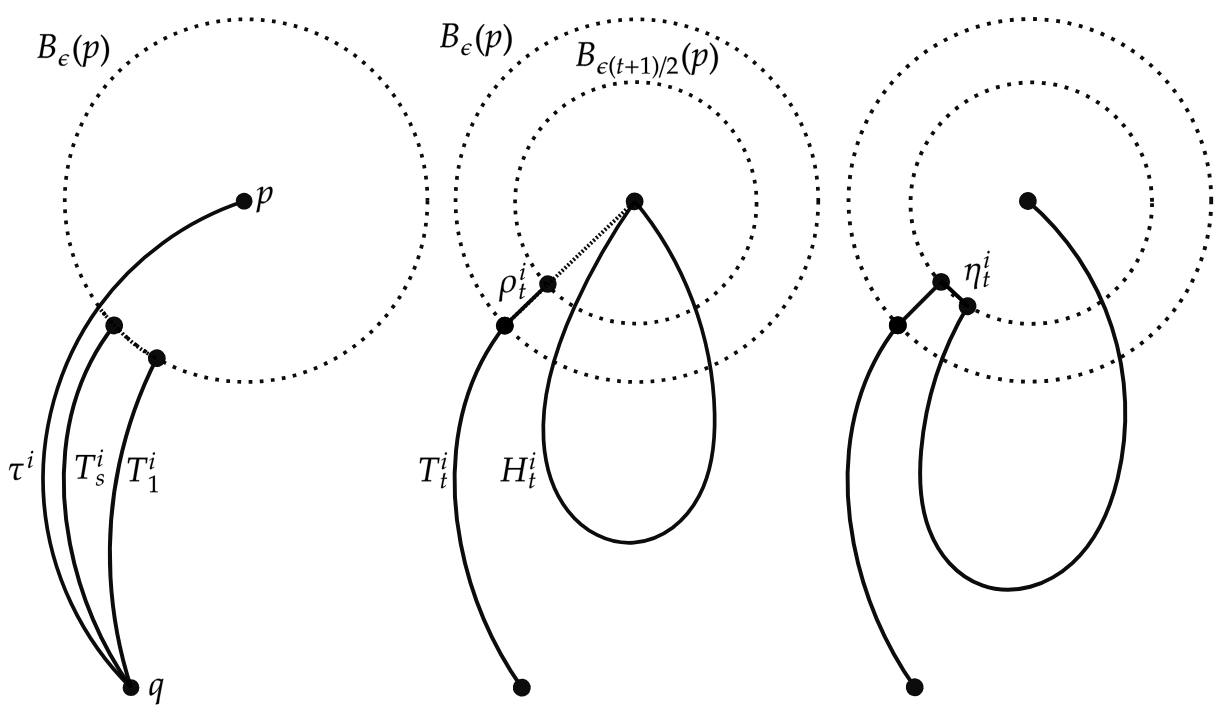}
    \caption{Proof of Lemma \ref{lemma:really_short_slicing}. Left: Defining $T_t^i$. Middle: Defining $\rho_t$. Right: Defining $\eta_t$ and a typical curve in the final sweepout.}
    \label{fig:example_hom}
\end{figure}
\subsection{Shortening Families}
\label{sec:families}
We now show how to apply the modified disk flow to one- or two-parameter families of $p$-admissible curves in order to obtain short simple geodesic loops at $p$. Because we want our two geodesic loops to have distinct images, we would like to identify any two curves that differ only in parameterization or orientation. However, the modified disk flow depends on the choice of which arc through $p$ contains the image of $t=0$. Recall that $\Pi_p M$ is the space of parameterized, oriented $p$-admissible curves in $M$. Let $\pi$ be the projection map that sends an oriented curve to the corresponding unoriented curve, remembering the choice of basepoint $t=0$, and define $\Sigma_pM=\pi(\Pi_p M)$.  
\par 
In order to prove our main theorem, we will need to apply the modified disk flow to representatives of homology classes of $\Sigma_p M$. Therefore consider a continuous map
$f:(X,\partial X)\to (\Sigma_p(M),\{p\})$, where $X$ is a $k$-dimensional finite simplicial complex and $\{p\}$ is the point curve at $p$. The properties of the image of $f$ under the modified disk flow, denoted $f_1$, satisfies the following properties which we will need to apply the Lyusternik--Schnirelmann argument in the subsequent section.
\begin{theorem}
\label{theorem:flow} Fix $l>0$.
    Suppose that every $p$-admissible curve of length at most $l+\delta$ is mapped to a curve of length at most $l-\delta$ by the modified disk flow.
    Let $f:(X,\partial X)\to (\Sigma_p M,\{p\})$ be a continuous map and suppose there exists $k>0$ such that every curve in the image of $f$ has at most $k$ vertices. If every curve in the image of $f$ has length at most $l+\delta/2$, then there exists a continuous map $f_1:(X,\partial X)\to (\Sigma_p M,\{p\})$ satisfying the following.
    \begin{enumerate}        
    \item 
        Every curve in the image of $f_1$ has length at most $l-\delta/2$.
        \item 
        The maps $f$ and $f_1$ are homotopic. 
        \item 
        There is some $k_1>0$ such that every curve in the image of $f_1$ has at most $k_1$ vertices. This constant depends only on $f$ and the choice of metric ball cover of $M$.
    \end{enumerate}
\end{theorem}
\noindent 
We will prove this theorem in the end of this section after we establish the definition of the modified disk flow on families.
\par
\begin{figure}[ht]
	\centering
	\includegraphics[width=\textwidth]{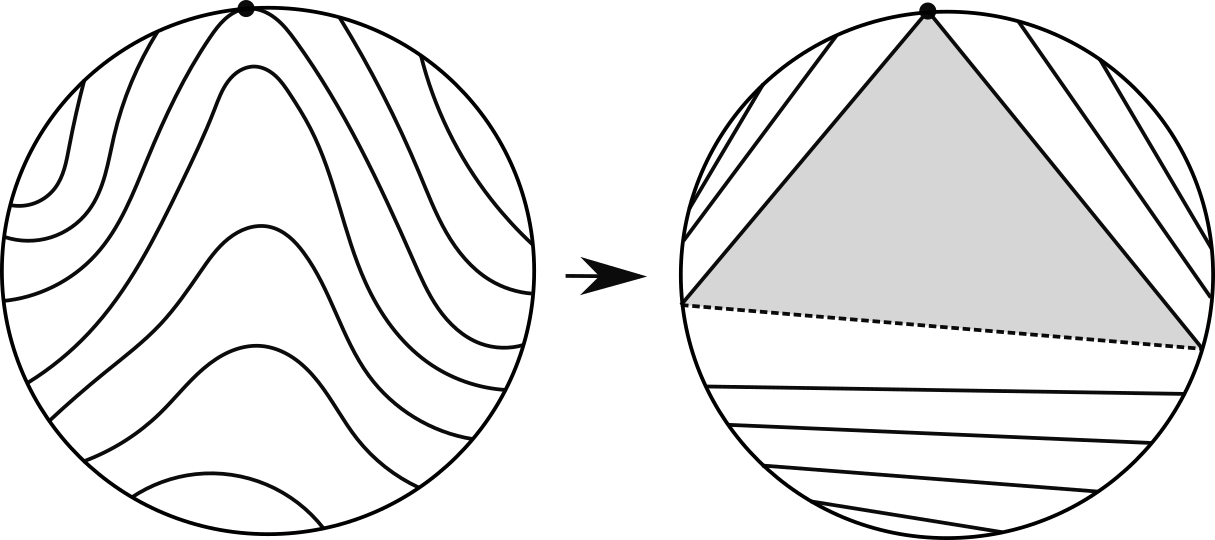}
	\caption{A family of curves in a disk before and after applying the (modified) disk flow. One of the curves is tangent to the boundary, creating a gap in the resulting family.}
	\label{fig:filling}
\end{figure}
\begin{figure}[ht]
    \centering
    \includegraphics[width=\textwidth]{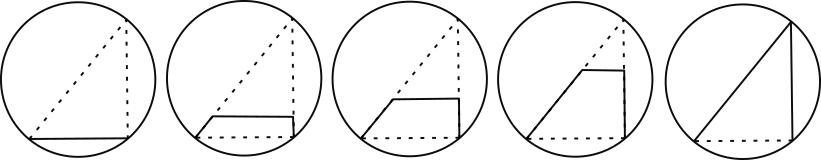}
    \caption{An example of an interpolating homotopy between edges of a geodesic triangle.}
    \label{fig:filling_1d}
\end{figure}

We cannot simply define the modified disk flow on families by applying the modified disk flow to each curve individually because of the existence of discontinuities. As in the original disk flow, if a curve in a family is ever tangent to one of the circles $\partial B_i$, there may be discontinuities in the disk flow process (see Figure \ref{fig:filling}). This is because, under the (modified) disk flow, an arc with a tangency will be replaced by a pair of geodesic segments, while a nearby arc that lies entirely within the disk will be replaced by only one geodesic segment. 
Thus there will be a discontinuity at any parameter where a curve has an interior tangency that ``moves inside" the metric ball. In order to fix this issue, we take advantage of the fact that we assumed our metric to be analytic in Lemma \ref{lemma:digons_and_loops} to prove that these discontinuities form a reasonably regular subset of the parameter space.
\begin{lemma}
    \label{lemma:analytic}
    Suppose that the metric on $M$ is analytic. Let $f:X\to \Sigma_p M$ be a continuous family of curves of length at most $L$. Suppose moreover that, for some integer $k>0$, every curve in the image of $f$ is a piecewise geodesic with at most $k$ vertices. For sufficiently small $\epsilon>0$, there exists a continuous family $g:X\to \Sigma_p M$ homotopic to $f$ that satisfies the following properties. 
    \begin{enumerate}
        \item 
        For any $x\in X$, every vertex of the piecewise geodesic $f(x)$ lies within $\epsilon$ of a vertex of $g(x)$. Conversely, every vertex of $g(x)$ is within $\epsilon$ of a point (not necessarily a vertex) of $f(x)$.
        \item 
        The length of a curve in the image of $g$ is at most $L+\epsilon$.
        \item 
        The set of parameters $x\in X$ where the curve $g(x)$ has an interior tangency to the boundary of a given metric ball $B_i$ is a triangulable subset of a $(\dim(X)-1)$-dimensional submanifold of $X$. 
    \end{enumerate}
\end{lemma}
\begin{proof}
    Any oriented piecewise geodesic, such as a curve in the image of $f$, is uniquely identified by the ordered list of its vertex points. This list is a point of $(x_1,x_2,\cdots,x_{k-1},x_k)\in M^k$, where the original curve is the piecewise geodesic connecting $p$ to $x_1$, then connecting $x_1$ to $x_2$, and so on up until we connect $x_k$ back to $p$. Since we are interested in unoriented curves, we will take the projection onto the quotient $M^k/\sim$ defined by identifying $(x_1,x_2,\cdots,x_{k-1},x_k)$ with $(x_k,x_{k-1},\cdots,x_{2},x_1)$. Then we can view $f$ as a map $f:X\to M^k/\sim$ given by the projection of $(f_1(x),\cdots,f_k(x))$. Note that because $M$ is analytic, this quotient space has an analytic structure.
    \par
    By Proposition 8 of \cite{grauert}, given any $\delta>0$ the map $f$ is homotopic to an analytic map $g:X\to M^k/\sim$ given by the projection of $(g_1(x),\cdots,g_k(x))$ satisfying 
    $$\sup_{x\in X} \sup_{1\leq i \leq k} d(f_i(x),g_i(x))<\delta.$$ 
    Therefore $g$ satisfies the first claimed property. Moreover, $g$ is homotopic to $f$ through unoriented piecewise geodesic curves curves satisfying the same distance bound. For example, we can define the geodesic in $M$ connecting $f_i(x)$ to $g_i(x)$ as $\rho_i(t)$ and consider the homotopy given at time $t$ by connecting $p$ to $\rho_1(t)$, then $\rho_1(t)$ to $\rho_2(t)$, and so on, until connecting $\rho_k(t)$ to $p$. Note that we have only established that these curves are piecewise geodesic, not that they are unoriented $p$-admissible curves. We will prove this later.
    \par 
    We can choose $\delta$ small enough that $g(x)$ is not longer than $L+\epsilon$. The vertex $g_i(x)$ is within $\delta$ of $f_i(x)$, so the geodesic arc connecting $g_i(x)$ and $g_{i+1}(x)$ has length no longer than the length of the arc between $f_i(x)$ and $f_{i+1}(x)$ plus $2\delta$. Recall that each $g(x)$ consists of $2k+2$ geodesic arcs. Therefore the length of $g(x)$ is at most $L+4(k+1)\delta$. Thus we can ensure the desired length bound holds by choosing $\delta<\epsilon/4(k+1)$.
    \par
    We now prove that the third condition holds. Define $\{x_i\}$ and $\{c_i\}$ such that the balls $B_i=B_{r_i}(c_i)$ form our chosen metric ball cover of $M$. Because the metric on $M$ is analytic, we can assume that each $r_i$ is small enough that the distance function $d(\cdot,c_i)$ is analytic on $B_i$. Consider the analytic function $h_{ij}:X\to\bR$ defined by $h_{ij}(x)=d(g_j(x),c_i)$. The set $h_{ij}^{-1}(r_i)$ corresponds to the set of curves in the family $g$ whose $j$th vertex lies on $\partial B_i$. By Sard's Theorem, we can slightly alter the radius of our metric ball cover to ensure that $r_i$ is a regular value of $h_{ij}$, and hence that $h^{-1}_{ij}(r_i)$ is a $(\dim(X)-1)$-dimensional submanifold of $X$. However, we want to consider only when this vertex intersects $\partial B_i$ as an interior tangency. Because we are considering a piecewise-geodesic curve, the point $g_j(x)$ is an interior tangency exactly when $g_{j-1}(x)$, $g_{j+1}(x)$, and $B_i$ are all on the same side of the line tangent to $\partial B_i$ at $g_j(x)$. Let $v$ be a vector tangent to $\partial B_j$ based at $g_j(x)$. Since $g_{j\pm1}(x)$ is a point on a surface, in the co-ordinates of $M$ its position relative to the tangent line at $g_j(x)$ is given by the sign of the determinant
    $
    \sigma_{j\pm 1}(x)$ of the matrix consisting of the vectors $v$ and $ g_j(x)-g_{j\pm1}(x) $.
    Because the manifold co-ordinates are analytic, this function is also analytic. Choose $v$ so that if $\sigma_{j\pm1}(x)>0$, then the corresponding point lies on the desired side of the tangent line. Therefore the set
    $$
    Z_{ij}=\{x\in X \mid h_{ij}(x)=r_i,\ \sigma_{j-1}<0, \text{ and } \sigma_{j+1}<0\}
    $$
    is the set of parameter values where has $g(x)$ an interior tangency to $\partial B_i$ at the point $g_j(x)$. Because this set is cut out by finitely many zero sets and inequalities of analytic functions defined on the compact set $X$, it is a subanalytic set. Therefore, $Z_{ij}$ admits a (locally finite) triangulation by \cite{Hironaka1974TriangulationsOA}, proving our claim.
    \par 
    Finally, we need to ensure that $g$ consists of $p$-admissible curves. We can do so by picking $\delta$ sufficiently small. Let $U_\delta^j(x)$ be the image of all minimizing geodesics that connect points in the $\delta$-neighbourhood of $g_j(x)$ to points in the $\delta$-neighbourhood of $g_{j+1}(x)$, where we define $g_0(x)=g_{k+1}(x)=p$. Recall that for every $x\in X$, $f(x)$ has no self-intersections except perhaps non-transversely at $p$. Therefore we can choose $\delta$ small enough to ensure that for every $x\in X$, $U_\delta^j(x)\cap U_\delta^i(x)=\emptyset$ unless $i=j\pm1$ or $f_j(x)=p$. Therefore $g(x)$ can self-intersect only along edges incident to $g_j(x)$ for those $j$ satisfying $f_j(x)=p$. If necessary, redefine $\delta$ so that $\delta<r_1$, where $r_1$ is the radius of the ball $B_1$ in our cover centred at $p$. Apply the modified disk flow in the disk centred at $p$ of radius $\delta$ to all curves $g(x)$ individually. Note that the modified disk flow does not depend on the orientation of a curve, so this is a well-defined operation. This disk is large enough to contain all self-intersections of all curves $g(x)$, but small enough that no curve $g(x)$ has an inner tangency to its boundary (and hence no discontinuities can occur). Therefore the resulting map is a continuous family of $p$-admissible curves of length at most $L+\epsilon$, and has an identical set of points tangent to the interior of $\partial B_i$ as the original map.
\end{proof}
Now that we know the set of discontinuities is reasonable, we will show how to remove the discontinuities. At a discontinuity caused by one tangency, the tangent curve is shortened as a pair of geodesics, while a nearby curve in the family that lies fully inside the ball is shortened as a single geodesic. Therefore we need to find a homotopy between these two curves that passes through short $p$-admissible curves (see Figure \ref{fig:filling_1d}). More generally, we need to be able to homotope between any edge of a geodesic polygon and the remaining edges. In fact, we want to be able to extend this homotopy to two-dimensional regions-- that is, we want to be able to remove discontinuities from 2-parameter families of curves defined on small polygons in parameter space. We can do these things by showing that the space of such curves is homotopy equivalent to a point. We explain how to do so in the following lemma.
\begin{lemma}
	\label{lemma:geo_polys_epsilon}
	Let $x_1,\cdots,x_k$ be a finite set of points on the boundary of a totally normal metric ball $B_r(x)\subset M$. Suppose $0<\epsilon$. Let $\Sigma_{r+\epsilon}$ be the set of simple piecewise geodesic curves $\gamma:[0,1]\to B_{r+\epsilon}(x)$ such that 
    \begin{enumerate}
        \item 
        $\gamma(0),\gamma(1)\in \partial B_r(x)$.
        \item 
        $\gamma(0)$ is within $\epsilon$ of $x_1$ and $\gamma(1)$ is within $\epsilon$ of $x_k$. 
    \end{enumerate}
    Let $\Sigma_{r+\epsilon}'\subset \Sigma_{r+\epsilon}$ be the subset of curves $\gamma$ such that 
    \begin{enumerate}
    \item
    Any intersection point $\gamma(t)\cap \partial B_r(x)$ is within distance $\epsilon$ of $x_i$ for some $1\leq i\leq k$.
    \item 
    Every vertex of $\gamma$ lies on either $\partial B_r(x)$ or $(\partial B_j\cap B_{r+\epsilon}(x))\setminus B_{r}(x)$ for some metric ball $B_j$ in our cover (possibly $B_j=B_r(x)$). 
    \end{enumerate}
    Then for sufficiently small $\epsilon$, there is a length non-increasing deformation retraction in $\Sigma_{r+\epsilon}$ of $\Sigma_{r+\epsilon}'$ onto the the space of geodesics in $B_r(x)$ connecting points within $\epsilon$ of $x_1$ to points within $\epsilon$ of $x_k$.
\end{lemma}
\begin{proof}
    We first prove that any edge of a convex geodesic polygon lying in $B_{r+\epsilon}$ can be homotoped with fixed endpoints to the remaining edges. Moreover, this homotopy depends continuously on the polygon and consists of simple curves that are not longer than the total length of the remaining edges and that lie within the polygon. Let $\rho_{a,b}$ denote the unique minimizing geodesic connecting any two points $a$ and $b$ in the polygon. 
	If we want our homotopy to end at the single edge $\rho_{x,y}$, let $u(t):[0,1]\to M$ be a parameterization of the remaining edges of the polygon in question. Then our homotopy is $t\mapsto u\mid_{[0,t]}*\rho_{u(t),u(1-t)}*u\mid_{[1-t,1]}$ (see Figure \ref{fig:filling_1d}). These curves are no longer than $u$, since we are replacing arcs of $u$ by minimizing geodesics. Moreover, they are simple and lie within the polygon by convexity. Note that also these homotopies are independent of the orientation of the underlying curve. 
    \par 
    We can generalize this statement to non-convex polygons in $B_{r+\epsilon}$. If the polygon has three vertices, it is necessarily convex. If it has four, then add a vertex between the two middle vertices and continue. Suppose the polygon has at least five vertices. In this case we apply the above homotopy simultaneously to the pair of (disjoint, convex) triangles defined by the first three vertices and last three vertices of the polygon. This reduces the number of total vertices of the polygon by two. We then repeat until we have reached zero vertices or reduced to a previous case. Note that, as above, this homotopy is independent of the orientation of the underlying curve, varies continuously with respect to the original curve, and consists of curves that are not longer than the remaining edges. However, these curves are not necessarily simple.
    \par
    We now prove the statement of the lemma. Let $\epsilon$ be small enough that the $\epsilon$-neighbourhood $U_i$ of $x_i$ contains the boundary of at most one other metric ball in our cover, say $B_{j_i}$. This is possible because the balls in our cover are in general position. Moreover, let $\epsilon$ be small enough that $B_{r+\epsilon}$ is totally normal.
    We start by defining a length non-increasing deformation retraction in $\Sigma_{r+\epsilon}$ of $\Sigma'_{r+\epsilon}$ onto $\Sigma'_{r}$. Consider a curve $\gamma \in\Sigma'_{r+\epsilon}$. An arc of $\gamma\cap U_i$ consists of arcs inside $B_r(x)$ and arcs outside of it. Because $\epsilon$ is small, the exterior arcs can only have vertices on $\partial B_j$ (recall that we consider $B_r(x)$ to be closed). Thus, to move $\gamma$ into $B_r(x)$ we need to be able to homotope any arc of $\gamma\cap U_i$ with endpoints on $\partial B_r(x)\cap U_i$ and vertices on $\partial B_j\cap U_i$ to the minimizing geodesic (in $B_r(x)$) connecting its endpoints. Moreover, these curves must be simple. We claim that the above arguments give the desired homotopy, where we view the arc and the minimizing geodesic connecting its endpoints as a geodesic polygon. Since this homotopy lies within the convex hull of the arc, it indeed passes through curves in $B_{r+\epsilon}$. It remains to see that the curves are simple in the case that the polygon is non-convex. It is enough to check that no subarc lies within the geodesic triangle whose vertices are the first three or final three vertices of the arc. Because the arc is simple and has all vertices on $\partial B_j$, this is indeed the case. 
    \par 
    We then once again apply our above homotopies in $\Sigma_r'$
    to produce a length non-increasing deformation retraction of $\Sigma_r'$ onto the space of geodesics in $B_r(x)$ connecting points within $\epsilon$ of $x_1$ and to points within $\epsilon$ of $x_k$. Since every curve in $\Sigma_r'$ has its vertices on $\partial B_r(x)$, our homotopies pass through simple curves even for arcs that form non-convex polygons.
\end{proof}
    \noindent An analogous version of this result holds for arcs with a fixed vertex at $p$. Note that this result extends naturally to unoriented curves because our homotopies do not depend on orientation. This is a fact which we will use in our proof of Theorem \ref{theorem:main} in the next section.
    \par 
    We now define the image of $f:(X,\partial X)\to (\Sigma_p M, \{p\})$ under the modified disk flow. We denote the image by $f_1$.
    Our first step is to approximate $f$, which we will temporarily view as a map $f:X\to\Sigma_p M$, using Lemma \ref{lemma:analytic}. Thus given any $\epsilon>0$ we have a map $g:X\to \Sigma_p M$ homotopic to $f$ such that the set of times where the $j$th vertex of the piecewise geodesic curve $g(x)$ has an interior tangency to $\partial B_{i}$ is contained in a triangulable subset, denoted $Z_{ij}$, of a $(\dim(X)-1)$-dimensional submanifold of $X$. Moreover, $g$ satisfies
    $$\max_{x\in X}L(g(x))\leq \max_{x\in X}L(f(x))+\epsilon.$$
    Our next step is to triangulate $X$. First, 
    triangulate the subanalytic set $\cup_{i,j} Z_{ij}$ and 
    extend this to a triangulation of $X$. Then refine this triangulation so that every curve corresponding to a parameter point in a given cell is within distance $\epsilon$ of all other curves in that cell in the sup norm. 
    We define the portion of $f_1$ in the ball $B_i$ on the 0-, 1- and 2-skeleton of our triangulation of $X$ as follows.
    \par
    \textbf{0-skeleton.}
    The map $f_1$ is defined on a point $x_0$ in the 0-skeleton by applying the modified disk flow in the ball $B_i$ to the curve $g(x_0)$. Note that this is well-defined since the outcome of the modified disk flow dos not depend on the orientation of the curve.
    \par
    \textbf{1-skeleton.}
    Consider an edge $e_t$, $t\in[0,1]$, in the 1-skeleton. For each $j$, $e_t$ is either disjoint from $Z_{ij}$, lies entirely within $Z_{ij}$, or intersects $Z_{ij}$ only at one or both of its vertices. In the second and third cases, the $j$th vertex of $g(e_t)$ is an interior tangency to $\partial B_i$ at some $t$ values. In the third case, it is possible that this tangency is moving ``inward" along the edge, creating a discontinuity in the disk flow. 
    Note that due to the construction of the 1-skeleton, the only place that tangencies can move inward is at vertices, not at interior points of an edge.
    \par 
    First we apply the disk flow in the ball $B_i$ to $g(e_t)$ to obtain $\tilde{g}(e_t)$ for each $t$. Consider a continuous family of arcs $\tilde{g}_t^j$ of $\tilde{g}(e_t)\cap B_i$ such that $\tilde{g}_0^j$ contains at least one vertex moving inward at $e_0$. 
    By the proof of Lemma \ref{lemma:geo_polys_epsilon}, there is a fixed-endpoint homotopy between $\tilde{g}_0^j$ and the minimizing geodesic connecting its endpoints through curves no longer than $\tilde{g}_0^j$. 
    There is a similar homotopy between $\tilde{g}_1^j$ and the geodesic connecting its endpoints. The endpoints of these geodesics are within distance $\epsilon$ of each other, so we can homotope between them by increasing the length by at most $2\epsilon$. Thus there is a homotopy $H^j_t$ between $\tilde{g}_0^j$ and $\tilde{g}_1^j$ through piecewise geodesics of length at most $\max_t{L(\tilde{g}_t^j)}+2\epsilon$ whose endpoints are within $\epsilon$ of the endpoints of $\tilde{g}_t^j$.
    We can move the endpoints of $H^j_t$ to align with those of $\tilde{g}_t^j$ by again increasing the length by at most $2\epsilon$. We then replace $\tilde{g}_t^j$ by the resulting family of arcs, which has length at most $\max_t{L(\tilde{g}_t^j)}+4\epsilon$. 
    \par 
    Repeat the above process for all arcs of $\tilde{g}(e_0)\cap B_i$ that have inward-moving tangencies. Then apply the analogous process for all arcs of $\tilde{g}(e_1)\cap B_i$ that have inward-moving tangencies whose discontinuities were not resolved at the previous step.
    For each arc we replace, we add $4\epsilon$ to the maximum length of $\tilde{g}(e_t)$. Therefore, since there are at most, say, $k$ tangencies, the maximum length of $f_1(e_t)$ is
    $L(\tilde{g}(e_t))+4k \epsilon.$ This completes our definition of $f_1(e_t)$, $t\in[0,1]$. Note that the curves $f_1(e_t)$ are $p$-admissible for all $t$. The reason that no self-intersections are introduced (except perhaps non-transversely at $p$) is because the original curves $g(e_t)$ were $p$-admissible, and hence there are no other arcs in the convex hull of any arc of $\tilde{g}(e_t)$ that does not pass through $p$. Therefore the homotopies of Lemma \ref{lemma:geo_polys_epsilon} do not cross any other arcs of $\tilde{g}(e_t)$.
    \par
    \textbf{2-skeleton.} 
    Consider a 2-simplex in the 2-skeleton of our triangulation.
    We claim that $f_1$ restricted to the boundary of this simplex (as defined in the previous step) is homotopic to a constant map and hence $f_1$ can be extended to the interior of the simplex. Consider a family of arcs in $B_i$ parameterized by a boundary edge of this 2-simplex. This family of arcs was constructed by applying the modified disk flow to a family of arcs in the image of $g$ and then possibly applying Lemma \ref{lemma:geo_polys_epsilon}. If Lemma \ref{lemma:geo_polys_epsilon} was not applied, then this family of arcs satisfies the hypotheses of Lemma \ref{lemma:geo_polys_epsilon}, and hence each arc can be homotoped to the geodesic connecting its endpoints as in the case of the 1-skeleton. This homotopy varies continuously along the edge. If Lemma \ref{lemma:geo_polys_epsilon} was applied, then this family is the image of a homotopy between an arc with a tangency and the minimizing geodesic connecting its endpoints. At each point on the simplex edge, we can continue the homotopy to move the arc associated to that point to the geodesic connecting the arc's endpoints. This homotopy varies continuously along the edge, since we are applying an identical homotopy from a continuously varying starting point. Combining these two cases, we have a continuous family of homotopies defined on the boundary of our 2-simplex. Each homotopy starts at an arc in $B_i$ and ends at the geodesic connecting its endpoints. The resultant family of geodesics can then be homotoped to the same geodesic, completing our homotopy between $f_1$ restricted to the boundary and a constant map. Note that the curves in this homotopy consist of simple piecewise geodesics that have a bounded number of breaks and bounded length. Therefore the same is true of the new curves in the interior of the simplex. Repeating this for each family of arcs completes the definition of $f_1$ on the 2-simplex.
\par
This completes our definition of $f_1$ in $B_i$. We define $f_1$ in $M$ by applying the above process in each ball of our cover. Note that the map $f_1$ is not unique, as it depends on our choice of $g$ and our triangulation of $X$.
\par 
We would like to consider $f_1:X\to \Sigma_p M$ as a map $f_1:(X,\partial X)\to (\Sigma_p M,\{p\})$, so we need to check that $f_1(\partial X)=\{p\}$. Point curves are not necessarily preserved under the modified disk flow-- for instance, point curves in the interior of a 1- or 2-cell may be removed if there are non-constant curves on the cell's boundary. However, we can fix this by adjusting $g$ and our triangulation of $X$. First, modify our triangulation of $X$ to be a common subdivision of our original triangulation and a triangulation of $\partial X$. Given $y\in \partial X$, the curve $g(y)$ lies in an $\epsilon$-neighbourhood of $p$ since $f(y)=p$. Therefore under the modified disk flow every individual curve $g(y)$ is going to be mapped to $\tilde{g}(y)=p$. Moreover, no such $g(y)$ is tangent to a metric ball in our disk. Therefore since $\partial X$ contains no 2-cells, $f_1(y)=\tilde{g}(y)$, so $f_1(y)=p$ for every $y\in \partial X$.
\par
We now prove our desired properties of the modified disk flow on families.
\begin{proof}[Proof of Theorem \ref{theorem:flow}]
    We first prove the claimed length bound. We know that given any $\epsilon>0$ we can define $f_1$ and $g$ and such that
    \begin{align*}
        \max_{x\in X} L(f_1(x))\leq \max_{x\in X} L(\tilde{g}(x))+c\epsilon 
    \end{align*} for some positive integer $c$, where $\tilde{g}(x)$ is the image of the individual curve $g(x)$ under the modified disk flow. Moreover,
    $$
    \max_{x\in X} L(\tilde{g}(x))\leq L(g(x))\leq L(f(x))+\epsilon\leq l+\delta/2+\epsilon.$$ 
    Taking $\epsilon$ small enough that $2c\epsilon<\delta$, we obtain $L(g(x))\leq l+\delta$ and hence $L(\tilde{g}(x))\leq l-\delta$, giving
    \begin{align*}
        \max_{x\in X} L(f_1(x))\leq l-\delta+\delta/2=l-\delta/2,
    \end{align*}
    as claimed.
    \par 
    We next want to define a homotopy $f_t:(X,\partial X)\to(\Sigma_p M,\{p\})$ between $f$ and $f_1$.
    As before, we can define it in each disk $B_i$ individually. First, suppose $i \not= 0$. For any curve $x\in X$, the arcs $f(x)\cap B_i$ are simple and disjoint. The arcs $f_1(x)\cap B_i$ are also simple and disjoint. Therefore we can view $f_1(x)\cap B_i$ as the image of $f(x)\cap B_i$ under a boundary-fixing diffeomorphism of the disk $B_i$. These diffeomorphisms depend continuously on $x$, since both $f$ and $f_1$ are continuous. Thus we have a continuous map of the parameter space $X$ into $\operatorname{Diff}(D, rel \partial)$, the space of boundary-fixing diffeomorphisms of the disk. Because $\operatorname{Diff}(D, rel \partial)$ is contractible, this map is homotopic to the constant map that sends every point of $X$ to the identity diffeomorphism. Call this homotopy $F : X\times[0, 1]\to \operatorname{Diff}(D, rel \partial)$. We thus obtain a homotopy from the image family $f_1$ in $B_i$ to the original family $f$ in $B_i$ by the map $f_t (x) \to (F (x, t))(f (x))$ for $x\in X, t \in [0, 1]$. 
    \par
    For $i=1$, we have to work a little harder to produce a diffeomorphism because the modified disk flow can map disjoint arcs to arcs that intersect non-transversely at $p$ and vice versa. Because our curves are $p$-admissible, we can continuously perturb any curve $f(x)$ and its image $f_1(x)$ inside $B_1$ so that they both consist of simple, disjoint arcs.
    However, we also want to ensure that our diffeomorphisms fix $p$ so that our intermediate maps pass through curves based at $p$. 
    After ensuring that our curves are simple, for each $f(x)$ there is exactly one arc of each curve through $p$. This arc splits $B_1$ into two disjoint disks, say $D^A(x)$ and $D^B(x)$. We can similarly define subdisks $D^A_1(x)$ and $D^B_1(x)$ with respect to $f_1(x)$. Moreover, we can continuously define these disks at each $x\in X$. For each $x$, the modified disk flow then restricts to a map $D^A(x)\to D_1^A(x)$ and a map $D^B(x)\to D_1^B(x)$. Each such map is a boundary-fixing diffeomorphism of the disk, and hence by our above argument we can obtain a homotopy $f_t^A$ between the set of arcs of $f(x)\cap D^A(x)$ to the set of arcs of $f_1(x)\cap D_1^A(x)$ and similarly for the other disks. To complete our homotopy, we finish by homotoping the arcs through $p$ of $f$ to the arcs through $p$ of $f_1$, for example by applying Lemma \ref{lemma:geo_polys_epsilon}. 
    \par
    Note that by Property 3 of Lemma \ref{lemma:modified_disk_flow}, there is a number $k_1$ (depending only on $f$ and our metric ball cover) such that each curve $f_1(x)$ has at most $k_1$ breaks.
\end{proof}

\subsection{Counting Loops}	
\label{sec:counting_loops}
We will now use the modified disk flow to prove the existence of two short simple geodesic loops at $p$ in analogy with the proofs of the Lusternik--Schnirelmann theorem given in \cite{hass1994, klingenberg_1978}.

\begin{proof}[Proof of Theorem \ref{theorem:main}]
	We are interested in chains in $\Pi_p M$ and cycles in $\Sigma_p M$ with $\bZ_2$ coefficients. If $M$ does not admit two simple geodesic loops at $p$ of length at most $2d$, then by Lemma \ref{lemma:really_short_slicing} we have a strictly monotone meridional slicing $\Gamma_t$, of $M$ through curves starting at $p$ and ending at some $q$. Moreover, $L(\Gamma_t)$ is at most $7d$ plus some arbitrarily small constant $c>0$. We will use $\Gamma_t$ to define two chains in $\Pi_p M$ which we call $u_0$ and $u_1$. Our cycles in $\Sigma_p M$ will be the images under $\pi$ of these chains.
	\par
	For convenience, we will extend $\Gamma_t$ periodically to all $t\in\bR$. Define $u_0:([0,1],\{0,1\})\to (\Pi_p M, \{p\})$ as follows. Fix some small $\epsilon>0$. For $t\in[\epsilon,1-\epsilon]$, take $u_0(t)=\Gamma_t*-\Gamma_{-t}$. For $t\in[0,\epsilon)$, let $u_0$ be a strictly monotone homotopy from $p$ to $\Gamma_\epsilon*-\Gamma_{-\epsilon}$ that increases lengths by at most some small $\delta$, where we take $\epsilon>0$ small enough that this is possible. Similarly, for $t\in(1-\epsilon,1]$, let $u_0$ be a strictly monotone homotopy from $\Gamma_{1-\epsilon}*-\Gamma_{-1+\epsilon}$ to $p$ that increases lengths by at most $\delta$. The second map $u_1:([0,1]\times[0,1/2],\{0,1\}) \to (\Pi_p M, \{p\})$ is essentially given by rotating $u_0$ through angle $\pi$. Specifically, for $(t,s)\in[\epsilon,1-\epsilon]\times[0,1/2]$, take $u_1(t,s)=\Gamma_{s+t}*-\Gamma_{-s-t}$. When $t\in[0,\epsilon)$, let $u_1$ be a strictly monotone homotopy from $p$ to $\Gamma_{s+\epsilon}*-\Gamma_{-s-\epsilon}$ that increases lengths by at most $\delta$, and when $t\in(1-\epsilon,1]$ let $u_1$ be a strictly monotone homotopy from $\Gamma_{s+1-\epsilon}*-\Gamma_{-s-1+\epsilon}$ to $p$. Moreover, we define both homotopies in a way that makes $u_1$ continuous. Finally, for notational convenience we reparameterize $u_0$ and $u_1$ to have domains $t\in[-1,1]$ and $(t,s)\in [-1,1]\times[0,1]$ respectively.
	\par 
	Define the $\bZ_2$-cycles $v_i=\pi(u_i)$. Note that this makes $u_1$ into a map $u_1:(X,\partial X)\to (\Sigma_p M,\{p\})$ where $X$ is the M\"obius band. We are concerned with the following two classes of cycles in $\Sigma_pM$. First, let $V_0$ be the smallest class of $\bZ_2$-cycles of the form $v:([-1,1],\{\pm1\})\to (\Sigma_pM, \{p\})$ that contains $v_0$ and contains all $\bZ_2$-cycles that are homotopic to $v_0$ as maps $([-1,1],\{\pm1\})\to (\Sigma_pM, \{p\})$. Similarly, $V_1$ is the smallest class of cycles $v:(X,\partial X)\to (\Sigma_pM, \{p\})$ that contains $v_1$ and contains all $\bZ_2$-cycles $v$ that are homotopic to $v_1$ as maps $(X,\partial X)\to (\Sigma_p M, \{p\})$.
	\par
	We want to find geodesic loops whose lengths realize the min-max values $l_i = \min_{v\in V_i}\max_\tau \{L(v(\tau))\}$. Note that $l_0\leq l_1$ because $u_0(t)=u_1(t,0)$ and hence any homotopy of $v_1$ determines a homotopy of $v_0$. Secondly, we claim that $l_0>0$ by the following argument. Let $f:S^2\to M$ be any map that is not nullhomotopic (e.g., the identity map). Suppose $u_0^\tau:([-1,1],\{\pm1\})\times[0,1]\to (\Sigma_pM, \{p\})$ is some homotopy with $u_0^0=u_0$. We want to parameterize $S^2$ by $u_0$ and then convert this homotopy to a homotopy $f^\tau$ of $f$. Every point in $S^2$ lies on some curve $u_0(t)$. First, perturb $u_0$ slightly near $q$ so that only one curve crosses $q$ and that all curves intersect only at $p$. Then for $x\not=p$ the choice of curve containing $x$ is unique, so if $x=(u_0(t))(s)$, we define $f^\tau(x)=u_0^\tau(t)(s)$. 
	Our map is also well-defined on $p$ via $f^\tau(p)=p$. 
	Now suppose that $l_0=0$, so that $u_0$ can be homotoped so that its image lies in an arbitrarily small disk around $p$. Then $u_0$ can be homotoped via some $u_0^\tau$ to the constant map $p$, and hence $f$ can be homotoped via some $f^\tau$ to the constant map $p$. This is a contradiction of the fact that $f$ is not nullhomotopic, and hence we see that $l_0>0$.
	\par 
	In addition, our definition of $u_0$ and $u_1$ shows that $l_i\leq 2\max_t\{L(\Gamma_t)\}=14d+\delta$ for any $\delta>0$, and hence $l_i\leq 14d$. In fact, $l_0\leq 8d$. This is because $V_1$ contains a cycle defined similarly to $v_0$ except using $\Gamma_0*-\Gamma_{-t}$ in place of $\Gamma_t*-\Gamma_t$, where without loss of generality $\Gamma_0$ is a shortest curve in the family $\Gamma_t$ (and hence of length at most $d$ by construction). 
	\par 
	Before we can prove anything else useful about $l_i$ and any corresponding critical points, we need the following intermediate lemma (c.f. Lemma 3.12 in \cite{hass1994}).
	\begin{lemma}
		\label{lemma:final_lemma}
		Fix $l>0$.
		Let $K_l$ be the (possibly empty) set of unions of simple geodesic loops at $p$ of total length $l$. Given any open neighbourhood $W_l$ of $K_l$ in $\Pi_p M$, there exists $\delta>0$ such that for any $p$-admissible curve $\gamma$ with $L(\gamma)<l+\delta$, either the modified disk flow maps $\gamma$ to a curve of length at most $l-\delta$ or $\gamma\in W_l$.
	\end{lemma}
	\begin{proof}
		Suppose the conclusion did not hold. Then we could find some $p$-admissible loops $\gamma_m$ so that
		\begin{enumerate}
			\item $L(\gamma_m)< l + 1/m$
			\item $L(\gamma_m')> l - 1/m$
			\item $\gamma_m\not\in W_l$
		\end{enumerate}
		where $\gamma_m'$ is obtained by applying the modified disk flow to $\gamma_m$ once. Because their lengths are uniformly bounded, there is a subsequence of the $\gamma_m$ that converges to some $\gamma_\infty$. As in the proof of Lemma \ref{lemma:modified_disk_flow}, the modified disk flow must fix the length of $\gamma_\infty$, or else for some large $m$ we would be able to contradict the assertion that $L(\gamma_m')> l - 1/m$. Therefore by Lemma \ref{lemma:modified_disk_flow}, $\gamma_\infty$ is either a point curve or a union of simple geodesic loops at $p$. Moreover, by the same argument we must have $L(\gamma_\infty)=\lim\limits_{m\to\infty}L(\gamma_m)=l$, since if $\gamma_\infty$ was much shorter than the $\gamma_m$ for large $m$, then $\gamma'_m$ is much shorter than $\gamma_m$.
		Note that this means that $\gamma_\infty$ cannot be a point curve, because in this case we could find some $\gamma_m$ of arbitrarily small length, and hence some $\gamma_m'$ of arbitrarily small length. Therefore $\gamma_\infty$ must be a union of simple geodesic loops, so since $L(\gamma_\infty)= l$, $\gamma_\infty$ is an element of $K_l$. In particular, $\gamma_\infty\in W_l$ and hence $\gamma_m\in W_l$ for sufficiently large $m$. This is a contradiction of our choice of $\gamma_m$. 
	\end{proof}
	Using this lemma, we will show that there exists a union of simple geodesic loops at $p$ with collective length realizing the min-max value $l_i$. Suppose there was no such set of loops. Define $\Sigma_p^lM$ as the set of elements of $\Sigma_pM$ of length at most $l$. Then by Lemma \ref{lemma:final_lemma}, there is some $\epsilon>0$ such that the modified disk flow maps $\Sigma_p^{l_i+\epsilon}$ into $\Sigma_p^{l_i-\epsilon}M$ (since we can take $W=\emptyset$). Consider a cycle $v\in V_i$ whose image consists of curves of length at most $l_i+\epsilon/2$. By Theorem \ref{theorem:flow}, there is a map $v'\in V_i$ such that the curves in the image of $v'$ have length at most $l_i-\epsilon/2$. This contradicts the definition of $l_i$, so the desired loops must exist.
	\par 
	Recall that we desire two distinct simple geodesic loops. If the union of loops realizing the value $l_0$ contains at least two distinct loops in its image, then there are two distinct simple geodesic loops at $p$ whose lengths sum to at most $l_0$. Therefore the shorter has length at most $l_0/2$, and the longer has length at most $l_0$. Otherwise, the critical value $l_0$ is realized by a single prime simple geodesic loop of length $l_0$. In this case we consider the union of loops realizing the value $l_1$. Again, if there are at least two distinct loops in its image, then there are two distinct simple geodesic loops at $p$ whose lengths are at most $l_1/2$ and $l_1$, although it is possible that one of these loops is the same as the one that realizes the critical value $l_0$. 
	The case of concern is when both $l_0$ and $l_1$ are realized by a single loop. If $l_0\not=l_1$, then these loops must be distinct (since they are prime) and we are done. Suppose instead $l=l_0=l_1$ and we obtain identical prime loops of length $l$. If there are any simple geodesic loops at $p$ of length strictly less than $l$, again we are done because this loop must be distinct from the (prime) loop of length $l$. If this is not the case, then we claim that there are in fact infinitely many distinct simple geodesic loops at $p$ of length $l$. We now prove this assertion.
	\par 
	Let $W_l\subset \Sigma_p M$ be an open neighbourhood of the set of unoriented simple geodesic loops at $p$ of length $l$. Note that this time we do not need to consider unions of geodesic loops, since we have assumed there are no geodesic loops of length less than $l$. Using Lemma \ref{lemma:final_lemma}, choose $\epsilon$ so that the modified disk flow sends $\Sigma^{l+2\epsilon}_pM$ into $W_l\cup\Sigma^{l-2\epsilon}_pM$. We then apply Theorem \ref{theorem:flow} to a cycle in $V_1$ of curves of length at most $l+\epsilon$, obtaining another cycle $v\in V_1$ of curves of length at most $l-\epsilon$. Pick a homotopically non-trivial
    curve $c(t):([0,1],\{0,1\})\to (X,\partial X)$.
    We claim that $\operatorname{img}c\cap v^{-1}(W_l)$ is non-empty. Because there are no simple geodesic loops at $p$ of length strictly less than $l$, the modified disk flow eventually maps every curve in $\Sigma^{l-\epsilon}_pM$ to the point curve $p$. Therefore if $\operatorname{img}c\cap v^{-1}(W_l)=\emptyset$, $v\circ c$ would map into $\Sigma^{l-\epsilon}_pM$ and hence we could use the modified disk flow to homotope $v\circ c$ to the constant map that sends every $x\in X$ to the point curve $p$. However, $v\circ c$ is homotopically non-trivial as it is homotopic to $v_0$. This is a contradiction, and hence $c\cap v^{-1}(W_l)\not=\emptyset$.
	\par
	Consequently, in fact $v^{-1}(W_l)$ supports a non-trivial $\bZ_2$ 1-cycle $k$, since otherwise we can assume $\partial (v^{-1}(W_l))$ consists of closed curves and we could concatenate arcs of $c\cap (v^{-1}(W_l))^c$ with arcs of $\partial (v^{-1}(W_l))$ to form a closed curve that does not intersect $v^{-1}(W_l)$. This is a contradiction unless this curve bounds a disk, in which case the original curve can be pushed entirely into $v^{-1}(W_l)$. If there are only finitely many simple geodesic loops at $p$ of length $l$, then we can choose $W_l$ small enough that each connected component contains only curves homotopic to a particular (unoriented) geodesic loop. Then $v\circ k$ lies entirely in one such connected component of $W_l$, and hence is nullhomotopic. This is a contradiction, as $v\circ k$ is homotopic to $v_0$. 
	\par 
	This proves that there must be two distinct simple geodesic loops at $p$ of length at most $\max\{l_0,l_1/2\}$ and $l_1$. Because $l_0 \leq 8d$ and $l_1\leq 14d$, this proves our main theorem.
\end{proof}

\section*{Acknowledgments}
The author would like to thank the anonymous referee for their many helpful suggestions. The author would also like to thank Alexander Nabutovsky and Regina Rotman for their invaluable guidance. This work was supported in part by an NSERC Canada Graduate Scholarships Doctoral grant. The author also gratefully acknowledges support from the Institute for Advanced Study Summer Collaborators 2024 program.
\bibliographystyle{alpha}
\bibliography{bibliography}
\end{document}